\documentclass[11pt,twoside]{article} 
\usepackage{amsmath}
\usepackage{latexsym} 
\usepackage{diagbox}
\usepackage{rotating}
\usepackage{enumitem}
\input{amssym.def} 
\input{amssym} 
\newfont{\fra}{eufm10 scaled 1095} 
\newfont{\Bb}{msbm10 scaled 1095} 
\newfont{\Bbg}{msbm10 scaled 1680} 
\newcommand\CC{{\mbox{\Bb C}}} 
\newcommand\RR{{\mbox{\Bb R}}} 
\newcommand\NN{{\mbox{\Bb N}}} 
\newcommand\ZZ{{\mbox{\Bb Z}}} 
\newcommand\A{{\mbox{\Bb A}}} 
\newcommand\Z{{\Bbb Z}}
 
\newcommand\QQ{{\mbox{\Bb Q}}}

\newcommand\fg{{\frak{g}}} 
\newcommand\fh{{\frak h}}

\newcommand\fm{{\frak m}} 
 
\newcommand\fk{{\frak k}}

\newcommand\bm{{{\bf m}}} 

\newcommand\cA{{\cal A}}
\newcommand\cB{{\cal B}}
\newcommand\cC{{{\cal C}}} 
 
\newcommand\cH{{{\cal H}}} 
\newcommand\cS{{\cal S}} 
\newcommand\cF{{\cal F}} 
 
\newcommand\cU{{\cal U}}
\newcommand\cM{{\cal M}}
\newcommand\ph{\varphi} 
 
\newcommand\eps{\varepsilon} 

\newcommand\fosc{\frak o \frak s \frak c} 
\newcommand\fspin{\frak s \frak p \frak i \frak n}

\newcommand{\Iso}{\mathop{{\rm Iso}}} 
\newcommand{\GL}{\mathop{{\rm GL}}} 
\newcommand{\SO}{\mathop{{\rm SO}}} 

\newcommand{\Spin}{\mathop{{\rm Spin}}} 
\newcommand{\Hom}{\mathop{{\rm Hom}}} 
 
\newcommand{\Osc}{{{\rm Osc}}} 
\newcommand{\cSpan}{{\overline{\rm span}}} %
\newcommand{\GG}{{\Gamma\backslash G}}
\newcommand{\Id}{{{\rm id}}} 
\newcommand{\ad}{{{\rm ad}}} 
 
\newcommand{\tr}{\mathop{{\rm tr}}} 
\newcommand{\Ad}{\mathop{{\rm Ad}}} 
\newcommand{\sgn}{\mathop{{\rm sgn}}}

\newcommand{\Ker}{\mathop{{\rm ker}}}

\newcommand{\sym}{\mathop{{\rm sym}}} 
 
\renewcommand{\Im}{\mathop{{\rm Im}}}

\newcommand{\spec}{{\rm spec}}
\newcommand{\Span}{{{\rm span}}} 
\newcommand\ip{{\langle\cdot \,,\cdot \rangle}} 
\newcommand\ipr{{(\cdot \,,\cdot )}}

\newcommand\la{{\langle}} 
\newcommand\ra{{\rangle}} 
\newcommand\proof{{\sl Proof. }} 
\newcommand{\qed}{\hspace*{\fill}\hbox{$\Box$}\vspace{2ex}} 
\newtheorem{theo}{Theorem}[section] 
\newtheorem{pr}[theo]{Proposition}
\newtheorem{de}[theo]{Definition}
\newtheorem{ex}[theo]{Example}
\newtheorem{re}[theo]{Remark}

\newtheorem{ft}[theo]{Fact}
\begin{document} 
\title{The cubic Dirac operator on compact quotients of the oscillator group} 
\author{Ines Kath and Margarita Kraus}
\maketitle 
\begin{abstract}
  
We determine the spectrum of Kostant's cubic Dirac operator $D^{1/3}$ on locally symmetric Lorentzian manifolds of the form $\Gamma\backslash \Osc_1$, where $\Osc_1$ is the four-dimensional oscillator group and $\Gamma\subset \Osc_1$ is a (cocompact) lattice. Moreover, we give an explicit decomposition of the regular representation of $\Osc_1$ on $L^2$-sections of the spinor bundle into irreducible subrepresentations and we determine the eigenspaces of $D^{1/3}$.  
\end{abstract}

{\bf MSC2020:} 53C50, 35Q41, 58J50, 22E27 \\[1ex]
{\bf Keywords:}  Lorentzian manifold, cubic Dirac operator, locally symmetric space 
\tableofcontents
\section{Introduction}
This paper is a contribution to the spectral theory of the cubic Dirac operator on compact locally symmetric Lorentzian manifolds. The cubic Dirac operator has been introduced by Kostant \cite{Ko}  on naturally reductive spaces as a purely algebraic object.
However, it can be also considered as a geometric Dirac operator $ D^{1/3}$ that belongs to a family of Dirac operators $D^t$. These operators are induced by non-standard connections \cite{A,Go}. 
The square of the cubic Dirac operator satisfies a simple formula, actually it equals minus the Casimir operator up to terms of order zero~\cite{Ko}.

The spectrum of geometric Dirac operators on pseudo-Riemannian manifolds has been calculated for  some specific  examples. For the Dirac operator $D = D^{1/2}$ associated with the Levi-Civita connection, the point spectrum of the
pseudo-Riemannian torus $T^{p,q}$ has been computed, see \cite{B} for $T^{1,2}$
 and \cite{Lan}
for the general case. Kunstmann \cite{Kunst} studied the spectrum for pseudo-Riemannian spheres. In the case of even dimension of the manifold or even index of the metric, he computed the point spectrum and proved that the imaginary axis belongs to the continuous spectrum and that the residual spectrum is empty.
Reincke \cite{Rei} explicitly computed the full spectrum
of $D$  on  $R^{p,q}$, the flat torus $T^{p,q}$ and products of the form $T^{1,1}\times F$, where $F$ is an
arbitrary compact, even-dimensional Riemannian spin manifold.

Here we will consider the Dirac operator on four-dimensional compact homogeneous spaces  $G/\Gamma$, where $G$ is a solvable Lie group endowed with a bi-invariant Lorentzian metric and $\Gamma $ is a (cocompact) lattice in $G$. More exactly, $G$ will be the four-dimensional oscillator group. 

Homogeneous spaces of the form $G/\Gamma$ for solvable $G$ are not only examples of compact homogeneous Lorentzian manifolds but play a central role in their classification.
Indeed, Baues and Globke \cite{BG} proved the following result.
Let $M=G/H$ be a compact homogeneous pseudo-Riemannian manifold, and let $G$ be connected and solvable. Then $H$ is a lattice in $G$ and the pseudo-Riemannian metric on $M$ pulls back to a bi-invariant metric on $G$. In the Lorentzian case such Lie groups $G$ admitting a biinvariant metric were  classified by Medina \cite{Me}.  They are products of an abelian group by a so-called oscillator group, which is a certain semi-direct product of a Heisenberg group by the real line. Combining these results of Baues, Globke and Medina, Revoy, one obtains a classification of all solvable Lie groups $G$ for which there exist compact Lorentzian $G$-homogeneous spaces. This classification in the Lorentzian case can already be found in \cite{Z}. 
If we restrict ourselves to four-dimensional manifolds, we now see that the group $G$ either is isomorphic to the abelian group $\RR^4$ or to the four-dimensional oscillator group $\Osc_1$. Thus $M$ is a flat Lorentzian torus or a quotient of $\Osc_1$ by a lattice and the metric on $M$ is induced by the bi-invariant metric on $\Osc_1$.

Let us explain the four-dimensional oscillator group $\Osc_1$ in more detail. 
This group is a semi-direct product of the 3-dimensional Heisenberg group $H$ by the real line $\RR$, where $\RR$ acts trivially on the centre $Z(H)$ of $H$ and by rotation on $H/Z(H)$. In particular, it is solvable. As mentioned above, it admits a bi-invariant Lorentzian metric. This metric is a particular case of a plane wave metric.  

The Lie algebra $\fosc_1$ of $\Osc_1$
is spanned by a basis $X,Y,Z,T$ where $Z$ spans the centre and the remaining basis elements  satisfy the relations $[X,Y]=Z$, $[T,X]=Y$, $[T,Y]=-X$. This Lie algebra is strongly related to the one-dimensional quantum harmonic oscillator. Actually, the Lie algebra spanned by the differential operators $P:=d/dx$, $Q:=x$, $H=(P^2+Q^2)/2$ and the identity $I$ is isomorphic to $\fosc_1$. 

The oscillator group contains lattices.  Each lattice $L$ in $\Osc_1$ gives rise to a compact locally-symmetric Lorentzian manifold $L\backslash \Osc_1$. Notice  that for a lattice, we take the quotient of the left action and therefore write  the subgroup on the left side.
The problem of classifying lattices in $\Osc_1$ was first considered by Medina and Revoy \cite{MR}. Note, however, that  the result in \cite{MR}  is not correct due to a wrong description of the automorphism group of an oscillator group. Lattices of $\Osc_1$ (as subgroups) were classified up to automorphisms of $\Osc_1$ by Fischer \cite{F} and up to inner automorphisms of $\Osc_1$ by Fischer and Kath \cite{FK}.
Here we will not consider arbitrary lattices but we will concentrate on basic lattices, see Section  \ref{S4} for a justification of this assumption.

Let $X$ be a quotient of $\Osc_1$ by a basic lattice $L$. On $X$, we fix a spin structure and consider the spinor bundle $\Sigma$. We use the canonical indefinite inner product on smooth sections of $\Sigma(X)$ to define a Krein space $L^2(\Sigma(X))$. Then the cubic Dirac operator is defined on $L^2(\Sigma(X))$ and $iD^{1/3}:L^2(\Sigma(X))\rightarrow L^2(\Sigma(X))$ is essentially self-adjoint (as an operator on a Krein space). Our aim is to determine the spectrum of this operator.

 The scalar curvature of $X=L\backslash\Osc_1$ is zero. This implies that the zero order terms vanish in the formula for the square of $D^{1/3}$. Thus $(D^{1/3})^2=-\Omega$, where $\Omega$ denotes the Casimir operator on $L^2(\Sigma(X))$. In order to determine the point spectrum of $D^{1/3}$, we will consider $L^2(\Sigma(X))$ as a representation of $\Osc_1$. We provide an explicit decomposition of 
$L^2(\Sigma(X))$ into a (discrete) direct sum of irreducible subrepresentations with finite multiplicities, see Theorem~\ref{theo}.
The cubic Dirac operator preserves each summand and its square $-\Omega$ acts by scalar multiplication on it. So we determine first the eigenvalues of $-\Omega$. To do so, we use results of \cite{FK}. We consider a finite cover $\tilde X$ of $X$ such that the pull back of the spin structure of $X$ to $\tilde X$ becomes trivial. Then we decompose $L^2(\Sigma(\tilde X))=L^2(\tilde X)\otimes \Delta$ according to \cite{FK} and afterwards we determine the space of sections that are invariant under all decktransformations of the cover $\tilde X\to X$.
This indeed allows to determine the point spectrum of $-\Omega$, see Theorem~\ref{spec}. In order to determine the point spectrum of $D^{1/3}$, we finally show that all square roots of eigenvalues of $-\Omega$ are eigenvalues of $D^{1/3}$. Using the decomposition of $L^2(\Sigma(X))$, we can also prove that the whole spectrum of $D^{1/3}$ is all of $\CC$. From general properties of Dirac operators, it follows, moreover, that  the residual spectrum is empty. In summary, we obtain:

\begin{theo}\label{T1}
Let $L\subset \Osc_1$ be a basic lattice and $X:=L\backslash \Osc_1$ be the quotient space. Let $D^{1/3}:L^2(\Sigma(X))\rightarrow L^2(\Sigma(X))$ be the cubic Dirac operator for a fixed spin structure on $X$. The spectrum of $D^{1/3}$ is equal to $\CC$. The point spectrum of $D^{1/3}$ depends on the spin structure. 
It  consists of 0 if the kernel of $\Omega$ is non-trivial and the two roots of each of the non-zero eigenvalues of $-\Omega$ (see Theorem~\ref{spec} for an explicit description of the point spectrum of $-\Omega$  depending on the spin structure). 
The residual spectrum is empty. 
\end{theo}

 In particular, it turns out that the point spectrum for a basic lattice is always discrete. 
This is no longer true for arbitrary lattices.
In Section \ref{S64}, we give examples of shifted lattices for which the point spectrum of $D^{1/3}$ on the quotient has accumulation points.

We not only determine the eigenvalues of $D^{1/3}$, but we also describe the corresponding eigenspaces explicitly, see Section~\ref{S62}.

In Section \ref{S65}, we use these results to determine the spectrum also for  all other $D^t$. We describe the point spectrum of $D^t$ in terms of the eigenvalues of $D^{1/3}$ and show that the full spectrum of $D^t$ is also equal to $\CC$ and that the residual spectrum of $D^t$ is empty. 

\section{Basic notions}\label{S2}
\subsection{Invariant connections on Lie groups}
In this short subsection, we want to recall some basic facts on connections on Lie groups. Let $G$ be a simply-connected Lie group endowed with a bi-invariant semi-Riemannian metric $\ip_G$. Let $\fg$ denote the Lie algebra of $G$. The metric $\ip_G$ on $G$ corresponds to an $\ad(\fg)$-invariant scalar product $\ip$ on $\fg$. As a Lie group, $G$ is endowed with a canonical left-invariant connection $\nabla^0$, for which all left-invariant vector fields are parallel. We use $\nabla^0$ to define a one-parameter family of metric connections $\nabla^t$ by $$\nabla^t_X Y=\nabla^0_X Y+ t[X,Y]=t[X,Y]$$ for left-invariant vector fields $X,Y$. For $t=1/2$, we obtain the Levi-Civita connection of $\ip_G$. All these connections define the same divergence operator on vector fields since $\ip_G$  is bi-invariant.

Let us explain the relation between the connection $\nabla^t$ of $G$ 
 and the canonical connections of $G$ understood as a reductive homogeneous space $K/H$ for $K:=G\times G$  and $H:=\Delta G\subset G\times G$. Let $\fk$ and $\fh$ denote the Lie algebras of $K$ and $H$, respectively. Then
\[\fm_t:=\{((t-1)v,tv)\mid v\in\fg\}\subset\fk,\ t\in\RR,\]
is a one-parameter family of $\Ad(H)$-invariant complements of $\fh$ in $\fk$. Each of the decompositions $\fk=\fh\oplus\fm_t$ defines a canonical connection on $G=K/H$ in the sense of Wang. 
This canonical connection coincides 
with the  left-invariant connection $\nabla^t$ of the Lie group $G$ introduced above. For more information, see \cite[Section 5.2]{book}.

\subsection{Spin structures on quotients of Lie groups}
Here we gather some facts on spin structures on quotients of Lie groups by discrete subgroups.
Let $G$ be a simply-connected Lie group endowed with a bi-invariant semi-Riemannian metric $\ip_G$, an orientation and a time orientation. As above, let $\fg$ denote the Lie algebra of $G$ and $\ip$ the induced scalar product on $\fg$. 

Let $\Gamma\subset G$ be a discrete subgroup and consider the quotient $\Gamma\backslash G$. The metric, the orientation and the time orientation on $G$ induce a metric, an orientation and a time orientation on the quotient. The tangent bundle of $\GG$ is trivial. Let $\pi:G\to\GG$ denote the projection, then
\[\GG \times \fg \longrightarrow T(\GG),\quad (\Gamma g, X)\longmapsto d \pi dL_g (X)\]
is an isomorphism of vector bundles.
Consequently, the bundle of oriented and time-oriented orthonormal frames on $\GG$ equals 
\[P_{\rm SO^+}(\Gamma\backslash G)=\GG\times {\rm SO}^+(\fg), \] 
where $\SO^+(\fg)$ denotes the identity component of $\SO(\fg):=\SO(\fg,\ip)$.
The spin structures on $\Gamma\backslash G$ are classified by $H^1(\Gamma\backslash G,\ZZ_2)$, see \cite{B}, Satz 2.6. or \cite{Fr}. Since $$H^1(\Gamma\backslash G,\ZZ_2) = \Hom(H_1(\Gamma\backslash G),\ZZ_2) = \Hom(\pi_1(\Gamma\backslash G),\ZZ_2) = \Hom(\Gamma,\ZZ_2),$$  
we obtain a one-to-one correspondence between isomorphism classes of spin structures on $\Gamma\backslash G$ and homomorphisms $\eps:\Gamma\rightarrow \ZZ_2=\ZZ/2\ZZ$, see \cite[Section 2.2]{Fr}  for the more  general case of a covering map.  The spin structure corresponding to $\eps$ equals 
$$P_{\Spin^+}^\eps(\Gamma\backslash G):=G\times_{\Gamma\hspace{-1pt},\,\eps } {\rm Spin}^+(\fg) :=\Gamma\backslash \big(G\times {\rm Spin}^+(\fg)\big),$$
where the action of $\Gamma$ on $G\times \Spin^+(\fg)$ is given by 
$$\gamma\cdot (g,\alpha)=(\gamma g,(-1)^{\eps(\gamma)}\alpha),$$
see \cite[Folg.~2.3]{B} or \cite[Prop.\,1.4.2]{G}. Let $\Delta_\fg$ denote the spinor representation of $\Spin^+(\fg)$. We obtain 
$$\Sigma(\Gamma\backslash G)=G \times_{\Gamma\hspace{-1pt},\,\eps } \Delta_\fg$$
for the spinor bundle associated with $P_{{\rm Spin}^+}^\eps(\Gamma\backslash G)$. We identify vector fields on $\GG$ with functions $X:G\to \fg$ that are invariant under left translation by $\Gamma$ and we
identify smooth sections of $\Sigma(\Gamma \backslash G)$ with  smooth functions $\psi:G\to \Delta_\fg$ that are invariant under the action of $\Gamma$, that is 
\begin{equation}\label{inv}
\psi(\gamma g)=(-1)^{\eps(\gamma)} \psi(g)
\end{equation}
for all $\gamma\in\Gamma$.

\subsection{Krein spaces}
Since the natural scalar product on the spinor bundle of a Lorentzian spin manifold is indefinite, sections of the spinor bundle do not constitute a Hilbert space in a natural way.  Therefore we will work in Krein spaces as it is done in \cite{B}.  For a general theory of such spaces see \cite{Bo, L}.

Let $K$ be a complex vector space and $\ip$ a possibly indefinite inner product on $K$. We define symmetric operators and selfadjoint operators  on $K$ in the same way as in the definite case. 

\begin{de} A Krein space $(K,\ip)$ consists of a complex vector space $K$ and an indefinite inner product $\ip$ on $K$ such that there exists a selfadjoint linear map $J:K\to K$ with the following properties: 
\begin{enumerate}
    \item $(\cdot \,,\cdot):=\langle \cdot\,,J\cdot\rangle$ is a positive definite inner product that makes $K$ a Hilbert space,
    \item  $J^2 = \Id$.
\end{enumerate} 
\end{de}
A linear map $J$ that satisfies this condition is called a {\it fundamental symmetry}. On $K$, we consider the strong topology. It is 
defined to be the norm topology of the Hilbert space $(K,(\cdot \,,\cdot))$, where $(\cdot \,,\cdot)=\langle \cdot\,,J\cdot\rangle$ for any fundamental symmetry $J$.   Although, in general, the linear map $J$ is not uniquely determined by $(K,\ip)$, the strong topology is well defined, i.e., independent of $J$.

If $A$ is a closed linear operator with a dense domain, then $\spec(A)$ denotes
the spectrum of $A$ and $\spec_p (A)$, $\spec_c(A)$, and $\spec_r(A)$ denote the discrete, continuous, and residual spectra, respectively.  

\begin{ft}{\rm \cite{L}}\label{fact}
If $A$ is a closed selfadjoint operator on a Krein space, then the complex conjugate of $\spec(A)$ satisfies 
\begin{enumerate}
    \item $\overline{\spec_r(A)}\subset \spec_p(A)$,
    \item $\overline{\spec_p(A)}\subset \spec_p(A)\cup \spec_r(A)$,
    \item $\overline{\spec_c(A)}\subset \spec_c(A)$.
\end{enumerate}
\end{ft}

\subsection{The space of spinors}\label{spinors}
In this section, we recall some basic facts on inner products on the space of spinor fields on   a pseudo-Riemannian manifold. 
Since later on we  will be interested in Dirac operators on Lorentzian manifolds, we restrict the explanations to the case where the metric of the manifold has Lorentzian signature. 
Moreover, we will concentrate on the case where the manifold is a quotient of a Lie group $G$ by a discrete subgroup $\Gamma$ although most of the results could be stated as well for general Lorentzian  manifolds.

Then, on the spinor module  $\Delta_\fg$, there exists a scalar product $\ip_\Delta$ of split signature such that 
\begin{equation}\label{invip}
 \la X\cdot u,v\ra_\Delta =\la u, X\cdot v\ra_\Delta 
\end{equation}
for all $X\in\fg$, where `$\cdot$' denotes the Clifford multiplication. This scalar product is unique up to multiplication by a real number different from zero. It defines a scalar product on the bundle $\Sigma(\Gamma\backslash G)$, which we also denote by $\ip_\Delta$. We choose a time-oriented left-invariant vector field $\xi$ on $G$ with $\la \xi,\xi\ra_G=-1$. This vector field defines a vector field on the quotient $\Gamma\backslash G$, which we also denote by $\xi$. We use $\xi$ to define a map $$J_\xi:\Sigma(\Gamma\backslash G)\longrightarrow \Sigma(\Gamma\backslash G),\quad u\longmapsto \xi\cdot u$$
and a positive definite scalar product $\ipr_{\Delta, \xi}:=\langle \cdot ,J_\xi ( \cdot) \rangle_\Delta$.  It satisfies 
$$(\xi \cdot u,v)_{\Delta,\xi}=(u,\xi\cdot v)_{\Delta,\xi}\,,\quad (X\cdot u,v)_{\Delta,\xi}=-(u,X\cdot v)_{\Delta,\xi}\,, \ X\in\xi^\perp. $$ 
The stabiliser of a timelike vector is  a maximal compact subgroup of the Lorentz group.  Therefore the vector field $\xi $ defines a reduction of the frame bundle of $\Gamma\backslash G$ to a maximal compact subgroup of the Lorentz group.
The scalar product  $\ipr_{\Delta, \xi}$ is invariant under this subgroup.

The volume form of the metric $\ip_G$ induces a measure $\mu$ on $\GG$, which is invariant under $G$. We define inner products on the space of compactly supported smooth sections of $\Sigma(\Gamma\backslash G)$ by 
$$\langle \ph,\psi\ra:=\int_{\Gamma\backslash G} \la \ph,\psi\ra_\Delta d\mu ,$$
$$( \ph,\psi)_\xi:=\int_{\Gamma\backslash G} ( \ph,\psi)_{\Delta, \xi} d\mu=\la \ph, J_\xi\circ\psi\ra.$$
The first one is indefinite, the second one is positive definite. 

We can identify spinors with smooth functions with values in $\Delta_\fg$ satisfying (\ref{inv}). If we identify, in addition,  $(\Delta_\fg,\ipr_{\Delta, \xi} )$ with the standard unitary space by choosing an orthonormal basis, then the scalar product $\ipr _\xi$ on smooth sections of $\Sigma(\Gamma\backslash G)$  becomes the standard $L^2$-product on functions (with several components). 

We define $L^2_\xi(\Sigma(\Gamma\backslash G)) $ as the completion of the space of compactly supported smooth sections in  $\Sigma(\Gamma\backslash G)$ with respect to the norm induced by $\ipr_\xi$. 
 
We want to compare the spaces $L^2_\xi(\Sigma(\Gamma\backslash G))$ for different choices of $\xi$. Let $r_\xi$ be the left-invariant Riemannian metric on $G$ defined by reversing the sign of $\ipr_G$ in direction of $\xi$. More exactly,
$$r_\xi(s\xi+X,t\xi+Y)=st+\la X,Y\ra_G$$ for $X,Y$ in the orthogonal complement of $\xi$ with respect to $\ip_G$. Let $\xi_1$ and $\xi_2$ be time-oriented left-invariant vector fields with $\la \xi_1,\xi_1\ra=\la\xi_2,\xi_2\ra=-1$. Since $r_{\xi_1}$ and $r_{\xi_2}$ are left-invariant, they are quasi-isometric, i.e., there exists a constant $C>0$ such that 
$$ \textstyle\frac1C r_{\xi_1}(X,X)\le r_{\xi_2}(X,X)\le Cr_{\xi_1}(X,X) $$ for all $X\in TG$. Of course, also the metrics induced by $r_{\xi_1}$ and $r_{\xi_2}$ on $\Gamma\backslash G$ are quasi-isometric. This implies that 
the spaces $L^2_{\xi_1}(\Sigma(\Gamma\backslash G))$ and $L^2_{\xi_2}(\Sigma(\Gamma\backslash G)) $ are the same in the following sense. They coincide as vector spaces (whose elements are equivalence classes of Cauchy series in the space of compactly supported smooth sections of $\Sigma(\Gamma\backslash G)$), and the identity $I : L^2_{\xi_1}(\Sigma(\Gamma\backslash G)) \to L^2_{\xi_2}(\Sigma(\Gamma\backslash G)) $ is a bounded isomorphism with bounded inverse \cite[Theorem 3.8]{Rei}.

Let us fix $\xi$ as above and put $K:=L^2(\Sigma(\Gamma\backslash G)):=L^2_{\xi}(\Sigma(\Gamma\backslash G))$ as a vector space. The map $J_\xi$ can be extended to $K$, and we can define an indefinite inner product $\ip=(\cdot\,,J_\xi\,\cdot)$. Then $(K,\ip)$ is a Krein space \cite[Satz 3.16]{B}. Its definition is independent of $\xi$. For any time-oriented left-invariant vector field $\xi'$ the map $J_{\xi'}$ is a fundamental symmetry.

\subsection{The cubic Dirac operator} \label{S25}
Every  connection $\nabla^t$ on $G$ induces a connection on $\Sigma(\Gamma\backslash G)$, which we will also denote by $\nabla^t$.
Consider a smooth section  $\psi:G\to \Delta_\fg$ of $\Sigma(\Gamma\backslash G)$, see (\ref{inv}). Let $X\in\fg$ be a left-invariant vector field on $G$. Then $X$ can also be considered as a vector field on $\Gamma\backslash G$.  Then $\nabla^t_X \psi = X(\psi) +t \gamma(X)\cdot \psi$, where
$$\gamma(X)=\frac 14 \sum_{a,b}\la X,[e_a,e_b]\ra e^a e^b$$
\cite[p.\,152]{M}.  Here and in the following, $\{ e_a\mid a=1,\dots,n\}$ denotes a basis of $\fg$ and $\{e^a\mid a= 1,\dots,n\}$ its dual basis with respect to $\ip$. These elements of $\fg$ can also be understood as vector fields on $\Gamma\backslash G$ or, equivalently, as 
constant maps $G\to\fg$, $g\mapsto e^a$ and $G\to\fg$, $g\mapsto e_a$. Then the Dirac operator corresponding to $ \nabla ^t$ is equal to  
\[D^t=\sum e^a\cdot \nabla^t_{e_a}.\]
If we apply this to a smooth section $\psi$ of the spinor bundle, we obtain
\begin{equation} 
D^t\psi = \sum e^a\cdot e_a(\psi)+  \frac{3t}{2} \sum_{a<b<c}\la[e_a,e_b], e_c\ra e^ae^be^c \cdot\psi,
\label{Diract} \end{equation}
where $e^ae^be^c$ is understood as an element of the Clifford algebra \cite[Eq.\,(5)]{A}.

For $t=1/3$ we obtain the cubic Dirac operator $D^{1/3}$. The square of this operator is related to the Casimir operator $\Omega=\sum e_a e^a\in\cU(\fg)$ with respect to $\ip$ by
$$(D^{1/3})^2=-\Omega -\frac1{24} \tr \ad(\Omega) =-\Omega +\frac1{24} \sum_{a,b}\la [e_a,e_b],[e^a,e^b]\ra=-\Omega +\frac16 {\rm Scal}.$$
Here ${\rm Scal}$ denotes the scalar curvature of $\ip _G$ on $\Gamma\backslash G$.

For the following remark, let us again concentrate on the Lorentzian case in order to simplify the exposition. If $\ip_G$ is Lorentzian, then $iD^t:L^2(\Sigma(\Gamma\backslash G))\to L^2(\Sigma(\Gamma\backslash G))$ is essentially selfadjoint in the Krein space $(L^2(\Sigma(\Gamma\backslash G)),\ip)$.  This can be seen as follows. We have noticed  that $\nabla^t$ defines the same divergence operator as the Levi-Civita connection $\nabla^{1/2}$. Furthermore, the Riemannian metric $r_\xi$ on $\Gamma\backslash G$ that is obtained by reversing the sign in direction of a time-oriented left-invariant vector field $\xi$ is complete. Indeed, $r_\xi$ is left-invariant on $G$, hence $(G,r_\xi)$ is a homogeneous Riemannian manifold and therefore complete. Hence $(\Gamma\backslash G, r_\xi)$ is also complete. Now the assertion follows from \cite[Satz~3.19]{B}. 
\subsection{The right regular representation}\label{S24}
Let $(G,\ip_G)$ be as above and let $\Gamma$ be a cocompact discrete subgroup of $G$.
The right regular representation $\rho$ of $G$ on $L^2(\Gamma\backslash G)$ is the unitary representation given by 
\begin{equation}\label{rrr}
(\rho(g)(\varphi))(x)=\varphi(xg).
\end{equation}
It is a classical result that $(\rho,L^2(\Gamma\backslash G))$ is a discrete direct sum of irreducible unitary representations of $G$ with finite multiplicities, see e.g. \cite{Wo07}.

Let $F$ be an automorphism of $G$. For a representation $(\sigma,V)$ of $G$ we define a representation 
\begin{equation}\label{pull}
F^*(\sigma,V):=(\sigma\circ F,V).
\end{equation}
Then 
\begin{equation}\label{pull2}L^2(\Gamma\backslash G)\stackrel{\sim}\longrightarrow F^* (L^2(F(\Gamma)\backslash G)),\quad f\longmapsto f\circ F^{-1}
\end{equation}
is an equivalence of representations.

Recall that a smooth section of the spinor bundle $\Sigma(\Gamma\backslash G)$ is identified with a smooth $\Gamma$-invariant function $\ph: G\to \Delta_\fg$. In this way we can also define an action of $G$ on $L^2(\Sigma(\Gamma\backslash G))$ by (\ref{rrr}). 
\section{The oscillator group and its Lie algebra}
\subsection{The oscillator group}\label{S31}
The 4-dimensional oscillator group is a semi-direct product of the 3-dimensional Heisenberg group $H$ by the real line. Usually, the Heisenberg group $H$ is defined as the set $H=\CC\times\RR$ with multiplication given by 
$$(\xi_1,z_1)\cdot(\xi_2,z_2)=(\xi_1+\xi_2,z_1+z_2+\textstyle{\frac12}\omega(\xi_1,\xi_2)),$$
where $\omega(\xi_1,\xi_2):=\Im(\overline{\xi_1}\xi_2)$. Hence in explicit terms, the oscillator group is understood as the set $\Osc_1=H\times\RR$ with multiplication defined by 
$$(\xi_1,z_1,t_1)\cdot(\xi_2,z_2,t_2)=(\xi_1+e^{it_1}\xi_2,z_1+z_2+\textstyle{\frac12}\omega(\xi_1,e^{it_1}\xi_2),t_1+t_2).$$

If we identify $\Osc_1\cong\RR^4$ as sets, then the Lebesgue measure is left- and right-invariant with respect to multiplication in $\Osc_1$.

Let us consider the automorphisms of this group.
For $\eta\in\CC$, let $C_\eta:\Osc_1\rightarrow \Osc_1$ be the conjugation by $(\eta,0,0)$. Then
$$C_\eta: (\xi, z, t)\longmapsto (\xi+\eta-e^{it}\eta,z+\textstyle \frac12\omega(\eta+\xi,\xi-e^{it}\eta),t).$$
Furthermore, we define an automorphism $T_u$ of $\Osc_1$ for $u\in\RR$ by
\begin{equation}\label{DTu}T_u: (\xi, z, t)\longmapsto (\xi,z+ut,t).
\end{equation}
Finally, consider an $\RR$-linear isomorphism $S$ of $\CC$ such that $S(i\xi)=\epsilon iS(\xi)$ for an element $\epsilon\in\{1,-1\}$ and for all $\xi\in\CC$. Then 
$\epsilon=\sgn(\det S)$ and also
\begin{equation}\label{ES}F_{S}: (\xi,z,t)\longmapsto (S\xi,\det(S) z,\epsilon t)
\end{equation}
is an automorphism of $\Osc_1$. Each automorphism $F$ of $\Osc_1$ is of the form 
$$F=T_u\circ C_{\eta}\circ F_{S}$$
for suitable $u\in\RR$, $\eta\in \CC$ and $S\in\GL(2,\RR)$ as considered above \cite{F}. Besides $C_\eta$ also $F_S$ is an inner automorphism if $S\in\SO(2,\RR)$.

In some of our computations we will use a slightly different multiplication rule for the oscillator group. It looks more complicated than the usual one but it will make the computations easier. We use the well known fact that the Heisenberg group $H$ is isomorphic to the set $H(1)$ of elements $M(x,y,z)$ parametrized by $x,y,z\in \RR$ with group multiplication 
\[M(x,y,z)M(x',y',z')=M(x+x',y+y',z+z'+xy').\]
We define an action ${l}$ of $\RR$ on $H(1)$ by
$$l(t)(M(x,y,z))=\textstyle M\big(x \cos t -y \sin t, x \sin t + y\cos t, z + \frac {xy}2 (\cos (2t) -1)+\frac {x^2-y^2}4\sin (2t)\big)$$
and consider the semi-direct product
\begin{equation}
\Osc_1^M:=H(1)\rtimes_l \RR. \label{GM}
\end{equation}
The image of an element $t\in\RR$ under the identification of $\RR$ with the second factor of $G$ in (\ref{GM}) is denoted by $(t)$.
It is easy to check that 
\begin{equation}
\phi:\Osc_1\rightarrow \Osc_1^M,\quad (x+iy,z,t)\mapsto M(-y,x,z-\textstyle{}\frac12 xy)(t)\label{Eiso}
\end{equation}
is an isomorphism.
Also here we can identify $\Osc_1^M\cong\RR^4$. Then $\phi $  preserves the Lebesgue measure.

\subsection{The oscillator algebra}
The Lie algebra $\fosc_1$ of the four-dimensional oscillator group is spanned by elements $Z,X,Y,T$, whose non-vanishing commutators are
$$[T,X]=Y,\ [T,Y]=-X,\ [X,Y]=Z.$$
The following result about the centre of the universal enveloping algebra of $\fg$ is known, see \cite{MuRi}. We give a short self-contained proof.

\begin{pr}\label{zentrum}
The centre $Z(\cU(\fg))$ of the universal enveloping algebra $\cU(\fg)$ of $\fg$ is generated by  $\Omega_0:=X^2+Y^2+2ZT$ and $Z\in\fg$.
\end{pr}
\proof Obviously, $\Omega_0$ belongs to $Z(\cU(\fg))$. By the Poincar\'e-Birkhoff-Witt Theorem, the symmetrisation map 
\begin{equation} \label{sym}
\sym: S(\fg)^\fg\rightarrow Z(\cU(\fg))
\end{equation}
induces an isomorphism $S(\fg)^\fg\cong {\rm gr}(Z(\cU(\fg)))$. Hence it suffices to show that $S(\fg)^\fg$ is generated by the centre of $\fg$ and the preimage of $\Omega_0$ under sym. Obviously, it suffices to show this for any $S^k(\fg)^{\fg}$.

Instead of $\fg$, we consider its complexification $\fg_{\Bbb C}$. The vectors
$Z,T$, $N_+:=X+iY$, $N_-:=X-iY$ constitute a basis of $\fg_{\Bbb C}$. Their non-vanishing Lie brackets are
$$[T,N_+]=-iN_+,\ [T,N_-]=iN_-,\ [N_+,N_-]=-2iZ.$$
In this new basis, we have $\Omega_0=N_+N_-+2ZT+iZ$. Hence
$\Omega_S:=N_+N_-+2ZT\in S(\fg)^\fg$
is a preimage of $\Omega_0$ under sym.
Let $\omega$ be in $S(\fg)^\fg$ and assume that $\omega$ is homogeneous. Then $\omega=\sum_{k,l}N_+^k p_{k,l}(Z,T)N_-^l$, where $p_{k,l}$ is a homogeneous polynomial in $Z$ and $T$. Since
$$\ad(T)(N_+^k p(Z,T)N_-^l)=i(l-k)N_+^k p(Z,T) N_-^l,$$
we obtain $\omega=\sum_{k=0}^n N_+^k p_k(Z,T)N_-^k$ if $\omega $ is invariant. Moreover, in this case
$$\ad(N_+)(N_+^k p(Z,T)N_-^k)=iN_+^{k+1}\frac\partial{\partial T} p(Z,T) N_-^k-2ikN_+^kZp(Z,T)N_-^{k-1} $$
yields $\frac{\partial}{\partial T}p_n(Z,T)=0$, thus $p_n(Z,T)=a_nZ^{m_n}$ for some $m_n\in\NN$. Now we consider $\omega':=\omega-a_n\Omega_S^n Z^{m_n}$. Then $\omega'$ is in $S(\fg)^\fg$ and of the form $\omega'=\sum_{k=0}^{n-1}N_+^k p'_k(Z,T)N_-^k$. We proceed inductively and obtain $\omega=a_n\Omega_S^nZ^{m_n}+a_{n-1}\Omega_S^{n-1}Z^{m_{n-1}}+\dots+a_0Z^{m_0}.$
\qed

\begin{re} {\rm The symmetrisation map (\ref{sym}) is not a homomorphism. For the sake of completeness let us determine Duflo's factor for $\fg$ although we will not use it in the present paper. See \cite{M} for a general introduction to this subject. For $\xi =zZ+n_+N_++n_-N_-+tT$, we have
$$\ad(\xi)=\left(\begin{array}{cccc} 
0&2in_-&-2in_+&0\\
0&-it&0&in_+\\
0&0&it&-in_-\\
0&0&0&0
\end{array}
\right)$$
with respect to the basis $Z, N_+, N_-, T$.  This gives
$J(\xi)=\det (j(\ad (\xi)))=j(-it)j(it),$
where $$j(z)=\frac{\sinh z/2}{z/2}.$$
Hence Duflo's factor equals
\begin{eqnarray*}
J^{1/2}(\xi)&=&\frac{\sinh(it/2)}{it/2}=\frac{\sin(t/2)}{t/2}.
\end{eqnarray*}
}\end{re}
\subsection{The biinvariant metric and the cubic Dirac operator}\label{Diracosc}
On $\Osc_1$, there exists a 2-parameter family of bi-invariant metrics. The metrics are defined by the ad-invariant scalar products on $\fosc_1$ given by 
$\Span\{X,Y\}\perp\Span\{Z,T\},$ and
$$\la X,X\ra=\la Y,Y\ra =r,\ \la X,Y\ra=\la Z,Z\ra=0,\ \la T,T\ra=s,\ \la Z,T\ra=r $$
for $r>0$ and $s\in\RR$. It is well known that there is only one bi-invariant Lorentzian metric on $\Osc_1$ up to isometric Lie group isomorphisms \cite{MR}. The above defined family of metrics arises as the orbit of such a metric under the action of the automorphism group of $\Osc_1$. The Casimir operator corresponding to the metric with parameters $r>0$ and $s\in\RR$ is equal to
 $\frac 1r(\Omega_0-s Z^2)$, where $\Omega_0 \in Z(\cU(\fg))$ is as defined in Proposition~\ref{zentrum}.

In the present paper, we consider the metric for $r=1$ and $s=0$, i.e.,
$$\la X,X\ra=\la Y,Y\ra=\la Z,T\ra=1,\ \la X,Y\ra=\la Z,Z\ra=\la T,T\ra=0.$$
The Casimir operator $\Omega$ of this metric equals
$$\Omega=\Omega_0=X^2+Y^2+2ZT.$$
The dual basis of $e_1=Z,\ e_2=X,\ e_3=Y,\ e_4=T$ is $e^1=T,\ e^2=X,\ e^3=Y,\ e^4=Z$. Hence
the scalar curvature of the metric  induced on $G$ vanishes since
$$\sum_{a,b}\la [e_a,e_b],[e^a,e^b]\ra=0.$$
This yields $(D^{1/3})^2=-\Omega$ for the square of the cubic Dirac operator.

Now we examine the actual Dirac operator  $D^{1/3}$. Since 
$ \sum_{a<b<c}\la[e_a,e_b], e_c\ra e^ae^be^c
= \la[X,Y],T\ra XYZ= XYZ,
$
we obtain from  (\ref{Diract}) 
\[
 D^{1/3}\psi=X\cdot X(\psi) +Y\cdot Y(\psi)+Z\cdot T(\psi)+T\cdot Z(\psi)+ {\textstyle\frac12}XYZ\cdot\psi.\]
Let $\Delta=\CC^4$ denote the spinor module of the metric Lie algebra $\fosc_1$. We can choose a basis $u_1,\dots,u_4$ of $\Delta$ such that the Clifford multiplication by $Z,X,Y$ and $T$ is given with respect to this basis by
\[
Z=\left(\begin{array}{cc} A&0 \\0 & A\end{array}\right),\quad T=\left(\begin{array}{cc} B&0 \\0 & B\end{array}\right),\quad X=\left(\begin{array}{cc} 0&C \\C &0 \end{array}\right),\quad Y=\left(\begin{array}{cc} 0&iC \\-iC &0 \end{array}\right),\]
where
$$A=\left(\begin{array}{cc} 0&0 \\ \sqrt 2&0 \end{array}\right),\quad B=\left(\begin{array}{cc}0 &-\sqrt 2 \\ 0&0 \end{array}\right),\quad C=\left(\begin{array}{cc} -i&0 \\0 & i\end{array}\right).$$
In particular,
\begin{equation}
  XYZ=i\left(\begin{array}{cc} A&0 \\0 & -A\end{array}\right).\label{xyz}
\end{equation}
According to section \ref{spinors}, there is an indefinite scalar product $\ip_\Delta$ on $\Delta$ satisfying (\ref{invip}) and this scalar product is uniquely defined up to a constant. We fix it by $\la u_1, u_2\ra_\Delta=\la u_3,u_4\ra_\Delta=1$, and $\la u_i,u_j\ra_\Delta=0$ for all other indices.

We choose the timelike left-invariant vector field $\xi=\frac1{\sqrt2}(Z-T)$ in order to define a definite $\fspin(3)$-invariant scalar product: $(u,v)_\Delta:=\la u, \xi \cdot v\ra_\Delta$. The vectors $u_1,\dots,u_4$ constitute an orthonormal basis with respect to $(\cdot\,,\cdot)_\Delta$.

\subsection{Unitary representations of the oscillator group}\label{S34}
 The irreducible unitary representations of $\fosc_1$ can be determined by applying a generalised version of Kirillov's orbit method. An explicit description of these representations can be found in \cite[\S 4.3]{Ki04}, where the oscillator Lie algebra is called diamond Lie algebra.
Let us recall this description. Note that the case $c<0$ in item (iii) does not appear in \cite{Ki04}. 
The infinite-dimensional representations will be given only on the Lie algebra level. 

Every irreducible unitary representation of the oscillator group is equivalent to one of the following representations, see also \cite{FK}:
\begin{enumerate}[label=(\roman*)]
\item $\cC_d:=(\sigma_{d},\CC)$, $\sigma_d (\xi,z,t)=e^{2\pi i d t}$, $d\in\RR$,
\item ${\cal S}_a^\tau:=(\sigma:=\sigma_a^\tau, L^2(S^1))$, $a>0$,  $\tau\in \RR/\ZZ\cong [0,1)$ where $\sigma$ is given by
\begin{eqnarray*}
\sigma_*(Z)(\ph)&=&0\\
\sigma_*(X+iY)(\ph)&=&2\pi i a e^{-it} \ph\\
\sigma_*(X-iY)(\ph)&=&2\pi i a e^{it} \ph\\
\sigma_*(T)(\ph)&=&\ph'+i\tau \ph 
\end{eqnarray*} 
for $\ph=\ph(t)\in C^{\infty}(S^1)\subset L^2(S^1)$. The orthonormal system $\phi_n:=e^{int}$, $n\in\ZZ$ satisfies
\begin{eqnarray*}
\sigma_*(X+iY)(\phi_n) &=&2\pi i a \phi_{n-1},\\
\sigma_*(X-iY)(\phi_n) &=&2\pi i a \phi_{n+1}\\
 \sigma_*(T)(\phi_n)&=&i(n+\tau)\phi_n.
 \end{eqnarray*} 
\item For $c>0$, $d\in\RR$, we consider the Hilbert space
$$\cF_c(\CC):=\Big\{\ph:\CC\to\CC \mbox{ holomorphic }\Big|\ \int_{\Bbb C} |\ph(\xi)|^2e^{-\pi c |\xi|^2} c\, d\xi <\infty\Big\}$$
with scalar product
\begin{equation}
\la \ph_1,\ph_2\ra=\int _{{\Bbb C}} \ph_1(\xi)\overline {\ph_2(\xi)} e^{-\pi c |\xi|^2} c\, d\xi\label{ipF}
\end{equation}
for $\ph_1,\ph_2\in \cF_c(\CC)$.
Then the representation $\sigma:=\sigma_{c,d}$ on $\cF_c(\CC)$ is given by
\begin{eqnarray*}
\sigma_*(Z)(\ph)&=&2\pi ic\ph\\
\sigma_*(X+iY)(\ph)&=&2\pi c\xi \ph\\
\sigma_*(X-iY)(\ph)&=&-2\partial \ph \\
\sigma_*(T)(\ph)&=&2\pi i d \ph -i\xi \partial\ph.
\end{eqnarray*}
The functions $\psi_n:=\frac{(\sqrt{\pi c} \xi)^n}{\sqrt{n!}}$, $n\ge0$, constitute a complete orthonormal system of $\cF_c(\CC)$ and we have
$$ \sigma_*(Z)(\psi_n)= 2\pi ic\psi_n,\ \sigma_*(T)(\psi_n)= (2\pi d-n)i\psi_n$$
and, for $A_+:=\sigma_*(X+iY)$ and $A_-:=\sigma_*(X-iY)$,
\begin{eqnarray}
&&A_+(\psi_n)=2\sqrt{\pi c(n+1)} \psi_{n+1},\ n\ge0, \label{A+}\\
&&A_-(\psi_0)=0,\ A_-(\psi_n)=-2\sqrt{\pi c n} \psi_{n-1}, \ n\ge1.\label{A-}
\end{eqnarray}
Furthermore, for $c<0$, $d\in\RR$, we consider 
$$\cF_{c}(\CC):=\Big\{\ph:\CC\to\CC \mbox{ anti-holomorphic }\Big|\ \int_{\Bbb C} |\ph(\xi)|^2e^{\pi c |\xi|^2} |c|\, d\xi <\infty\Big\}$$
with scalar product given by (\ref{ipF}) with $c$ replaced by $-c$,  now for $\ph_1,\ph_2\in \cF_{c}(\CC)$.
The representation 
$\sigma:=\sigma_{c,d}$ on $\cF_{c}(\CC)$ is given by
\begin{eqnarray*}
\sigma_*(Z)(\ph)&=&2\pi ic\ph\\
\sigma_*(X+iY)(\ph)&=&-2\bar \partial \ph\\
\sigma_*(X-iY)(\ph)&=& -2\pi c\bar\xi\ph\\
\sigma_*(T)(\ph)&=&2\pi i d \ph + i\bar\xi \bar\partial\ph.
\end{eqnarray*}
Here, the functions $\psi_n:=\frac{(\sqrt{\pi |c|} \bar\xi)^n}{\sqrt{n!}}$, $n\ge0$, constitute a complete orthonormal system and we have
$$ \sigma_*(Z)(\psi_n)= 2\pi ic\psi_n,\ \sigma_*(T)(\psi_n)= (2\pi d+n)i\psi_n.$$
Now, Equations (\ref{A+}) and (\ref{A-}) hold for $A_+:=\sigma_*(X-iY)$ and $A_-:=\sigma_*(X+iY)$.
We will use the notation 
$\cF_{c,d}:=(\sigma_{c,d},\cF_c(\CC))$ for all $c\not=0$ and $d\in\RR$.
\end{enumerate}
Let $F$ be an automorphism of $\Osc_1$. In (\ref{pull}), we defined the pullback of a representation of $G$ by $F$.
The following table shows (the equivalence class of) $F^*(\sigma,V)$ for the case that $V$ is one of the irreducible unitary representations of $\fosc_1$ and $F$ is one of the (outer) automorphisms $T_u$ or $F_S$ introduced in Section \ref{S31} by (\ref{DTu}) and (\ref{ES}). Note that $C_\eta$ as an inner automorphism does not change $V$. 
\begin{equation} \label{Fiso}
\mbox{   
\begin{small}
\renewcommand{\arraystretch}{2}
\begin{tabular}{|c||c|c|c|}
  \hline
  $F$&$F^*\cC_d$ & $F^*{\cal S}^\tau_a$ &$F^*{\cal F}_{c,d}$\\
  \hline \hline
  $T_u$ &$\cC_d$ &${\cal S}^\tau_a $&${\cal F}_{c,d+uc}$ \\
  \hline
  $F_S$ &$\cC_{\epsilon  d}$&${\cal S}^{\epsilon\tau}_{|\hspace{-1pt}\det S|^{1/2}a}$&${\cal F}_{\hspace{-1pt}\det(S) c,\epsilon d}$\\
  \hline
\end{tabular}
\end{small}
}
\end{equation}
\section{Straight and basic lattices}\label{S4}
In Section~\ref{S6} we will study the spectrum of the cubic Dirac operator on compact quotients of the oscillator group. More exactly, we consider quotients of $\Osc_1$ by  discrete uniform subgroups of $\Osc_1$. We will call such subgroups lattices. This is justified by the fact that the group $\Osc_1$ is solvable  and therefore a quotient by a discrete subgroup is of finite measure (for the measure inherited from Haar measure on $\Osc_1$) if and only if it is compact.  The lattices of the oscillator group are known. They were classified up to automorphisms of $\Osc_1$ by Fischer \cite{F}. Since here we are interested in the spectrum of the quotient and therefore in the right regular representation, we need a classification up to inner automorphisms, which can be found in \cite{FK}. 

To avoid too much technical effort, we will concentrate on straight lattices, where a lattice in $\Osc_1$ is called \emph{straight} if it is generated by a lattice in $H$ and an element $\delta$ of the centre of $\Osc_1$. It can be shown that each lattice in $\Osc_1$ contains a sublattice of finite index which is a straight, see \cite[Section 8]{F}. In other words, each lattice in $\Osc_1$ is virtually straight.

Moreover, we will assume that the lattice is unshifted and normalised in the sense of \cite{FK}. A lattice $L$ in $\Osc_1$ is called {\em normalised} if the projection of $L\cap H$ to $H/Z(H)\cong \RR^2$ has covolume one with respect to the standard metric of $\RR^2$. A normalised straight lattice is called {\em unshifted} if $\delta$ can be chosen in the $\RR$-factor of $\Osc_1=H\rtimes\RR$, i.e., $\delta=(0,0,2\pi\kappa)$.
 This leads us to the following definition. 
\begin{de}
A lattice of $\Osc_1$ is called a basic lattice if it is normalised and generated by a lattice in the Heisenberg group and an element $(0,0,2\pi\kappa)\in \RR\subset H\rtimes\RR$.
\end{de}
The additional assumptions to be normalised and unshifted are justified by the fact that each  straight  lattice can be normalised and shifted by (outer) automorphisms of $\Osc_1$. More exactly, the following holds. Let $\cM^{\mbox{\rm{\footnotesize strt}}}$ denote the set of all isomorphism classes of straight lattices of $\Osc_1$ with respect to inner automorphisms of $\Osc_1$ and let ${\cal B}\subset \cM^{\mbox{\rm{\footnotesize strt}}}$ be the set of isomorphism classes of basic lattices. For a basic lattice $L$ we define numbers $\kappa=\kappa(L)\in\NN_{>0}$ by $L=\la L\cap H, \delta=(0,0,2\pi\kappa)\ra$ and $r=r(L)\in\NN_{>0}$ by $L\cap Z(H)=r\ZZ$.
\begin{pr}\label{umr}
The map
\begin{eqnarray*}
{\cal B}\times\RR_{>0}\times \RR/\ZZ &\longrightarrow& \cM^{\mbox{\rm{\footnotesize strt}}} \\
(L,a,s)&\longmapsto& F_S(T_u(L)), \qquad u= \frac s{2\pi \kappa r},\ S=a\cdot I_2
\end{eqnarray*}
I is a bijection. Here $I_2$ denotes the identity on $\RR^2$. 
\end{pr}
\proof The assertion follows from Thm.~4.12 in \cite{FK}. Indeed, the property to be straight is invariant under automorphisms. Therefore we can restrict the bijections in \cite[Thm.~4.12 ]{FK} to straight lattices. In the notation of \cite{FK}, we thus obtain bijections from  $(\cM_0\cap \cM^{\mbox{\rm{\footnotesize strt}}})\times \RR/\ZZ=\cB\times\RR/\ZZ$ to $\cM_1\cap \cM^{\mbox{\rm{\footnotesize strt}}}$ and from $(\cM_1\cap \cM^{\mbox{\rm{\footnotesize strt}}})\times\RR_{>0}$ to $\cM^{\mbox{\rm{\footnotesize strt}}}$. It remains to check that the composition of these bijections has the form asserted in the proposition. Let $L$ be a basic lattice. Then we have  $q=1$ and $x_\delta=y_\delta=0$ in item 1 in \cite[Def.~4.8]{FK}, which implies $v=w=0$ and therefore $s_0=1$. The assertion follows.  \qed

The computation of the spectrum relies on the decomposition of the right regular representation into irreducible subrepresentations.  Once this decomposition is known for basic lattices, the decomposition for arbitrary lattices can be derived using Proposition~\ref{umr}. Indeed,
according   to (\ref{Fiso}) the decomposition of  $L^2(F_S(T_u(L))\backslash \Osc_1)$ can be computed from that of $L^2(L\backslash \Osc_1)$.  That is why  we focus on basic lattices here.
\begin{pr}{\rm \cite[Remark 5.3]{FK}} \label{sl} A basic lattice is isomorphic via an inner automorphism of $\Osc_1$ to a lattice generated by 
\begin{eqnarray*}
&&\textstyle l_1:=\big(\frac{1}{\sqrt\nu},0,0\big),\quad l_2:=\big(-\frac{\mu}{\sqrt\nu}+i\sqrt\nu,0,0\big),\quad l_3:=\big(0,\frac{1}{r},0\big),\quad\textstyle l_4:=\big(0,0,2\pi\kappa \big),
\end{eqnarray*}
for some $\mu,\nu\in \RR, \nu> 0$ and  $r,\kappa\in \NN_{>0}$. 
\end{pr}

We will denote this lattice by $L_r(\kappa,\mu,\nu)$. In \cite{FK}, it is denoted by $L_r(2\pi\kappa,\mu,\nu,0,0)$, but here we do not need the last two parameters since we only consider straight lattices. 
\begin{re}{\rm
The lattices $L_r(\kappa,\mu,\nu)$ and $L_{r'}(\kappa',\mu',\nu')$ differ by an inner automorphism of $\Osc_1$ if and only if $r=r'$, $\kappa=\kappa'$ and $(\mu,\nu)$ and $(\mu',\nu')$ are on the same orbit of the ${\rm SL}(2,\ZZ)$-action on the Poincar\'e half plane \cite[Thm.~4.15]{FK}. 
}\end{re}

Now let $L$ be a basic lattice.  By Prop.~\ref{sl} we may assume that $L=L_r(\kappa,\mu,\nu)$ in the following. Let us rewrite the first two generators using the matrix
$$ T_{\mu,\nu}:= \left(\begin{array}{cc} \sqrt\nu&\frac{\mu}{\sqrt\nu}\\0&\frac 1{\sqrt\nu} \end{array}\right).$$
Identifying $\CC\cong\RR^2$ and using the standard basis $e_1,e_2$ of $\RR^2$, we obtain
$$l_1=(T_{\mu,\nu}^{-1}e_1,0,0), \quad l_2=(T_{\mu,\nu}^{-1}e_2,0,0).$$

As an abstract group, the lattice $L$ is  isomorphic to the direct product of a discrete Heisenberg group $H_1^r(\ZZ):=\la l_1,l_2,l_3\mid l_1l_3=l_3l_1, l_3l_2=l_2l_3, l_1l_2l_1^{-1}l_2^{-1}=l_3^r \ra$ and $\ZZ$. 

We fix a spin structure on $X=L\backslash \Osc_1$. As explained in Section~\ref{S2}, it is determined by a homomorphism $\eps: L\rightarrow (\ZZ_2,+)$.  We will use the notation \[\eps_i:=\eps(l_i),\ i=1,\dots,4\] and write $\eps=(\eps_1,\dots,\eps_4)$. Note that a map $\eps: L\rightarrow \ZZ_2$ is a homomorphism if and only if $r\eps_3=0$. Let again $\Delta=\CC^4$ denote the spinor module of the metric Lie algebra $\fosc_1$.

\begin{re}\label{isometrie}{\rm If $\eps_1$ and $\eps_2$ are different, we may assume that $\eps_1=0$ and $\eps_2=1$ by changing $\mu$ and $\nu$ if necessary. Indeed, let be given the  lattice $L=L_r(\kappa,\mu,\nu)$ and a spin structure defined by $\eps=(\eps_1,\dots,\eps_4)$. We define $\mu'$ and $\nu'$ by $\mu'+i\nu'=-(\mu+i\nu)^{-1}$ and consider $L':=L_r(\kappa,\mu',\nu')$. In particular, we have 
$$l_1=(T_{\mu,\nu}^{-1}e_1,0,0),\ l_2=(T_{\mu,\nu}^{-1}e_2,0,0),\ l'_1=(T_{\mu',\nu'}^{-1}e_1,0,0),\ l'_2=(T_{\mu',\nu'}^{-1}e_2,0,0).$$
Let $S\in\SO(2)$ be the multiplication by $(\mu-i\nu)/|\mu-i\nu|$ on $\RR^2\cong\CC$. Then 
\[ T_{\mu',\nu'}^{-1}=ST_{\mu,\nu}^{-1} \cdot \mbox{\small$\begin{pmatrix}0&-1\\1&0 \end{pmatrix}$}.\]
Let $F_S:\Osc_1\to\Osc_1$ denote the inner automorphism defined by $S$ according to (\ref{ES}). Then $F_S(l_2)=l_1'$, $F_S(l_1)=(l_2')^{-1}$ and $F_S(l_j)=l_j'$ for $j=3,4$.
In particular, $F_S$ induces an isometry from $L\backslash \Osc_1$ to $L'\backslash \Osc_1$. If we pull back the spin structure on $L\backslash \Osc_1$ defined by $\eps$ by the inverse of this isometry, we obtain the spin structure on $L'\backslash \Osc_1$ defined by $\eps'=(\eps_2,\eps_1,\eps_3,\eps_4)$. In particular, the spectra of $D^{1/3}$ on $L\backslash \Osc_1$ with respect to $\eps$ and on  $L'\backslash \Osc_1$ with respect to $\eps'$ coincide.
}\end{re}
\section{The right regular representation for basic lattices}
\subsection{Strategy}\label{S51}
Let $L$ be a basic lattice of $\Osc_1$. We consider $X:=L\backslash \Osc_1$.  The aim of this subsection is to decompose the representation $L^2(\Sigma(X))$ of $\Osc_1$ into irreducible components. We want to apply the results of \cite{FK} for the decomposition of the right regular representation on $L^2$-functions. In \cite{FK}, the decomposition of $L^2(L\backslash \Osc_1)$ is determined for arbitrary lattices, where first the computation is reduced to the case of unshifted and normalised lattices and then explicit formulas are given in this case. In particular, \cite[Prop.~7.2]{FK} describes the decomposition for basic lattices. In order to apply these results, we consider a finite covering of $X$ such that the lifted spin structure becomes trivial. 

More exactly, 
we consider the covering $\tilde X=L'\backslash \Osc_1$ of $X$, where $L'$ is the subgroup of $L$ generated by $l_1^2,\dots, l_4^2$. Since $L$ is a basic lattice, this subgroup is normal. We obtain that $X=I  \backslash \tilde X$, where $I\subset\Iso(\tilde X)$ is the finite group generated by the actions of $l_1, \dots, l_4$ on $\tilde X$. The spin structure on $X$ lifts to a spin structure on $\tilde X$, which is now the trivial one since $2\eps_j=0$, $j=1,\dots,4$.  Therefore the associated spinor bundle of $\tilde X$ equals $\tilde X\times \Delta$ and  sections in this bundle can be identified with functions from $\Osc_1$ to $\Delta$ that are invariant under left translation by elements of the lattice $L'$. To recover the  sections in the spinor bundle of $X$ from these sections we have to find those sections in $\tilde X\times \Delta$ that are invariant under  the action of the group $I$ of decktransformations,  where this action is defined as follows.   Let $\ph=f\otimes u\in C^\infty(\tilde X)\otimes \Delta$ be a (local) section. 
Then $l. \ph=l^*f \otimes (-1)^{\eps(l)} u=(-1)^{\eps(l)} l^*f \otimes u$, where $l^*$ denotes the left translation by $l\in L$. 
Thus we can identify \[L^2(\Sigma(X))=L^2(\Sigma(\tilde X))^I=(L^2(\tilde X)\otimes \Delta)^I=L^2(\tilde X)^I\otimes \Delta,\]
where an element  $[l]\in I$ induced by $l\in L$ acts on $L^2(\tilde X)$ by 
\begin{equation}\label{act} [l]. f=(-1)^{\eps(l)} l^*f. 
\end{equation}
Consequently, we can obtain a decomposition of $L^2(\Sigma(X))$ into irreducible subspaces in the following way.  First we decompose $L^2(\tilde X)$ according to \cite{FK}. Then, for each isotypic component, we determine the subspace of sections that are invariant under the action of $l_1,\dots, l_4$ by (\ref{act}). Finally, tensoring by $\Delta$ gives the result. 

As said above, the explicit formulas in \cite{FK} for the decomposition work under the assumption that the lattice is normalised. However, note, that our new lattice $L'$ generated by
\begin{eqnarray*}
&&\textstyle l_1^2:=(2T_{\mu,\nu}^{-1}e_1,0,0), \quad l_2^2:=(2T_{\mu,\nu}^{-1}e_2,0,0),\quad l_3^2:=\big(0,\frac{2}{r},0\big),\quad\textstyle l_4^2:=\big(0,0,4\pi\kappa \big),
\end{eqnarray*}
is not normalised. Indeed, the projection of $L'\cap H$ to $H/Z(H)\cong \RR^2$ is generated by the projections of $l_1^2$ and $l_2^2$. Thus it has covolume 4. Therefore we apply the automorphism $F_S$ for $S=\frac12\cdot\Id$ to $\Osc_1$. We have $F_S( \xi, z, t)=(\frac12\xi,\frac14 z, t)$.   Then
$
F_S(l_1^2)=l_1, F_S(l_2^2)=l_2, F_S(l_3^2)=(0,\frac1{2r},0), F_S(l_4^2)=(0,0,4\pi\kappa)=l_4^2.
$
Hence $F_S(L')$ is the basic lattice $L_{r'}(\kappa',\mu,\nu)$ for $r'=2r$ and $\kappa'=2\kappa$. Thus the formulas in \cite{FK} apply to $F_S(L')$.

\begin{pr} \label{umrech} We have \begin{equation}\label{equr} L^2(\Sigma(X))\cong F_S^* (L^2(F_S(L')\backslash \Osc_1)^{F_S(I)})\otimes \Delta,\end{equation}
where $F_S^*$ is understood as the pullback of a representation as defined by (\ref{pull}) and $F_S(I)$ is the finite group of isometries of $F_S(L')\backslash \Osc_1$ generated by $[F_S(l_1)],\dots,[F_S(l_4)]$, where $[F_S(l_j)]$ acts on $L^2(F_S(L')\backslash \Osc_1)$ by $[F_s(l_j)]. f=(-1)^{\eps_j} F(l_j)^*f$ for $j=1,\dots,4$.
\end{pr}
\proof We have seen above that
$L^2(\Sigma(X))=L^2(\tilde X)^I\otimes\Delta=L^2(L'\backslash\Osc_1)^I\otimes\Delta$. Now the assertion follows from (\ref{pull2}) and the fact that the equivalence $f\mapsto f\circ F_S^{-1}$ maps $L^2(L'\backslash\Osc_1)^I$  to $L^2(F_S(L')\backslash \Osc_1)^{F_S(I)}$. 
\qed

\subsection{The decomposition of $L^2(\Sigma(X))$}

As explained in  Subsection \ref{Diracosc} we have to describe the decomposition  of $L^2(\Sigma(X))$ into irreducible subrepresentations up to equivalence. We will use the representations introduced in Subsection  \ref{S34}. 
To formulate the result we need the following notations.

Let
$$\|(k,l)\|_{\mu,\nu}:=\big(\nu k^2 +\frac1\nu(-\mu k+l)^2\big)^{\frac12}=\|T_{\mu,\nu}^{-1}(l,k)^\top\|,$$
and
 \begin{eqnarray*}
 \alpha(\mu,\nu,a)&:=&\# \{(k,l)\in\ZZ^2\mid \|(k,l)\|_{\mu,\nu}=  a\}\\
\alpha_0(\mu,\nu,a)&:=&\# \{(k,l)\in\ZZ^2\mid \|(k,l)\|_{\mu,\nu}=  a,  \mbox{ $k$ even, $l$ even }\},\\
&=&\alpha(\mu,\nu,a/2)\\
\alpha_1(\mu,\nu,a)&:=&\# \{(k,l)\in\ZZ^2\mid \|(k,l)\|_{\mu,\nu}=  a,  \mbox{ $k$ even, $l$ odd}\},\\
\alpha_2(\mu,\nu,a)&:=&\# \{(k,l)\in\ZZ^2\mid \|(k,l)\|_{\mu,\nu}=  a,  \mbox{ $k$ odd, $l$ odd }\}.
\end{eqnarray*}
Moreover,
$$\A(\mu,\nu):=\{ a\in\RR_{>0}\mid \alpha(\mu,\nu,a)\not=0\}. $$
\begin{theo}\label{theo} On $X=L_r(\kappa,\mu,\nu)\backslash \Osc_1$ we consider the spin structure given by $\eps=(\eps_1,\dots,\eps_4)$, where we assume $(\eps_1,\eps_2)=(0,1)$ if $\eps_1\not=\eps_2$. Then we have
$L^2(\Sigma(X))= 4\cH_0\oplus 4\cH_1$ with
\begin{eqnarray}
\cH_0&\cong& \bigoplus_{n\in {\Bbb Z}}\bm _{\cC}(n) \cC_{\frac n{4\pi \kappa}} \oplus \bigoplus_{a\in{\Bbb A}(\mu,\nu)}
 \bigoplus_{K=0}^{2\kappa-1} \bm_{\cS}(a,K)  {\cal S}_{a/2}^{K/(2\kappa)},\label{H0}\\
\cH_1&\cong&\bigoplus_{m\in \Z_{\not=0}} \bigoplus_{n\in\Z}\bm_\cF(m,n) \cF_{\frac {rm}2,\,\frac n{4\pi\kappa}}, \label{H1}
 \end{eqnarray}
 where 
\begin{eqnarray}
\bm_\cC(n)&=&\label{mcn} \left\{\begin{array}{ll} 1, & \mbox{if } \eps = (0,0,0,n)\in \ZZ_2^4, \\
0, & \mbox{else,}
\end{array}\right.\\
\label{msa}
\bm_{\cS}(a,K)&=& \left\{\begin{array}{ll}\alpha_0(\mu,\nu,a), &\mbox{if } \eps=(0,0,0,K)\in \ZZ_2^4, \\
\alpha_1(\mu,\nu,a), &\mbox{if } \eps=(0,1,0,K)\in \ZZ_2^4,\\
\alpha_2(\mu,\nu,a), &\mbox{if } \eps=(1,1,0,K)\in \ZZ_2^4,\\
0, & \mbox{else,}
\end{array}\right.
\\\label{mfmn}
\bm_\cF(m,n)&=&\left\{\begin{array}{ll} \frac{r|m|}2, & \mbox{if } (\eps_3,\eps_4)=(m,n)\in \ZZ_2^2, \\
0, & \mbox{else.}
\end{array}\right.
\end{eqnarray}
\end{theo}
\proof We proceed according to the strategy outlined in Subsection~\ref{S51}. By Proposition~\ref{umrech}, we have to determine the decomposition of $F_S^* (L^2(F_S(L')\backslash \Osc_1)^{F_S(I)})$, where $F_S( \xi, z, t)=(\frac12\xi,\frac14 z, t)$ and $L'$ is the lattice generated by $l_1^2,\dots,l_4^2$. It turns out that calculations on $\Osc_1^M$ are easier than on $\Osc_1$. Therefore we transform the lattice $F_S(L)$ by the isomorphism $\phi:\Osc_1\rightarrow \Osc_1^M$ defined in (\ref{Eiso}). We obtain
\begin{eqnarray*}
&\textstyle \gamma_1:=\phi\circ F_S(l_1)=M(0,\frac1 {2\sqrt\nu},0),\quad    \gamma_2:=\phi\circ F_S(l_2)=M(-\frac{\sqrt\nu}2,-\frac{\mu}{2\sqrt\nu},\frac{1}{8}\mu), &\\
&\textstyle \gamma_3:= \phi\circ F_S(l_3)=M(0,0,\frac{1}{4r})=M(0,0,\frac{1}{2r'}),\quad \gamma_4:= \phi\circ F_S(l_4)= (2\pi\kappa)=(\pi\kappa'),&
\end{eqnarray*}
where $r'=2r$ and $\kappa'=2\kappa$.
We denote the lattices $(\phi\circ F_S)(L)$ and $(\phi\circ F_S)(L')$ by $\Gamma$ and $\Gamma'$, respectively.

 The push-forward of a representation  
 $(\sigma,V)$  of $\Osc_1$ is a representation  $(\phi^{-1})^*(\sigma,V)=(\sigma\circ\phi^{-1},V)$  of $\Osc_1^M$. 
 In the following, we will identify these representations with each other and omit $(\phi^{-1})^*$ in the notation. In particular, we identify the representation 
$ L^2(F_S(L')\backslash \Osc_1)^{F_S(I)}= L^2(F_S(L')\backslash \Osc_1)^{F_S(L)}$ with  $L^2(\Gamma'\backslash \Osc_1^M)^\Gamma$.

It is  natural to use the push-forwards of the irreducible representations of $\Osc_1$ as models for the irreducible representations of $\Osc_1^M$. 
 Then the irreducible unitary representations of $\Osc_1^M$ are $(\phi^{-1})^*\cC_d$, $(\phi^{-1})^*{\cal S}^\tau_a$ and $(\phi^{-1})^*{\cal F}_{c,d}$. 
According to the above remark we simply write $\cC_d$, ${\cal S}^\tau_a$ and ${\cal F}_{c,d}$ instead of $(\phi^{-1})^*\cC_d$, $(\phi^{-1})^*{\cal S}^\tau_a$ and $(\phi^{-1})^*{\cal F}_{c,d}$. Our task now is to describe explicitly the irreducible subrepresentations of the right regular representation $L^2(\Gamma'\backslash \Osc_1^M)$ of $\Osc_1^M$. Here we can use the results of \cite{FK}. The representation $L^2(\Gamma'\backslash \Osc_1^M)$
is equivalent to $\cH_0'\oplus \cH_1'$, where the subrepresentations $\cH_0'$ and $\cH_1'$ are given as follows:
$$
\cH_0'\cong \bigoplus_{n\in {\Bbb Z}} \cC_{\frac n{2\pi \kappa'}} \oplus \bigoplus_{a\in{\Bbb A}(\mu,\nu)}
 \bigoplus_{K=0}^{\kappa'-1} \alpha(\mu,\nu,a)  {\cal S}_a^{K/\kappa'},$$
$$
\cH_1'\cong\bigoplus_{m\in \Z_{\not=0}} |m|r'\cdot \bigoplus_{n\in\Z} \cF_{r'm,\,\frac n{2\pi\kappa'}}. \label{PH1}
$$
The subspaces of $\cH_0'$ and $\cH_1'$ corresponding to the irreducible subrepresentations of $L^2(\Gamma'\backslash \Osc_1^M)$ in the above formulas are explicitly given as follows. For a function $f:\Osc_1^M\to \CC$, we will denote by $f(x,y,z,t)$ the image of $M(x,y,z)(t)$ under $f$. Then the representation $\cH_0'$ is the direct sum of subspaces 
$$ \cSpan\{ \ph_n(x,y,z,t):=\ph_n(t):=e^{2\pi i d t}\} \cong \cC_{d} $$
for $d=\frac n{2\pi\kappa'}$
and
\begin{eqnarray*}
\bigoplus_{\substack{k,l\in{\Bbb N}\\ \|(k,l)\|_{\mu,\nu}=a}}\cSpan\left\{ \phi_{l,n}^k(x,y,z,t)=\exp\big(2\pi i (x,y)T_{\mu,\nu}^{-1}(l,k)^\top\big) \ph_n(t)\ \big|\  n\in K+\kappa'\ZZ\right\}&&\\[-3ex]
\cong\alpha(\mu,\nu,a){\cal S}_{a}^{K/\kappa'} \hspace{3cm}&&
\end{eqnarray*}
for $a\in{\Bbb A}(\mu,\nu)$ and $K=0,\dots,\kappa'-1$.

Next we want to describe the decomposition of $\cH_1'$.  For $m\in\ZZ_{\not=0}$ and $n\in\ZZ$, we consider the subspace
$${\cal W}_{m,n}:=\cSpan\{\theta_{m,n,k} \mid k\in\ZZ_{r'|m|}=\ZZ/r'|m|\ZZ\}$$
for 
\begin{eqnarray*}
\theta_{m,n,k}(x,y,z,t)&:=& e^{2\pi ir'mz}e^{int/\kappa'}\sum_{j\in \frac  k{r'm}+\Z}
e^{-\pi r' |m|(x+j\sqrt\nu)^2}e^{\pi i r'jm(j\mu+2\sqrt\nu y)}.\\
\end{eqnarray*}

We denote by $A_+$ the ladder operator $X+\sgn(m)iY$. Then $\cH_1'$ decomposes into the direct sum of subspaces
$${\cal W}_{m,n}\oplus A_+({\cal W}_{m,n})\oplus A_+^2({\cal W}_{m,n})\oplus\dots\ \cong |m|r'\cdot  \cF_{r'm,\,\frac{n}{2\pi\kappa'}}$$
for $m\in\ZZ_{\not=0}$ and $n\in\ZZ$.

In order to obtain $L^2(\Gamma'\backslash \Osc_1^M)^\Gamma$, we determine the elements in the isotypic components of \[L^2(\Gamma'\backslash \Osc_1^M)\cong \cH'_0\oplus \cH'_1\] that are invariant under $ \gamma_1,\dots,\gamma_4$. 
We compute
\begin{eqnarray*}
\gamma_1\cdot M(x,y,z)(t)&=&\textstyle  M(x,y+\frac1{2\sqrt\nu}, z)(t), \\
\gamma_2\cdot M(x,y,z)(t)&=& \textstyle M(x-\frac{\sqrt\nu}2,y-\frac\mu{2\sqrt\nu }, z+\frac 18\mu-\frac{\sqrt\nu}2 y)(t),\\
\textstyle \gamma_3\cdot M(x,y,z)(t) &=&\textstyle  M(x,y,z+\frac1{2r'})(t),\\
\gamma_4\cdot M(x,y,z)(t)&=&\textstyle  M(x,y,z)(t+\pi\kappa'),
\end{eqnarray*}
 where we used that $\kappa'$ is even. 
Thus, the action defined by (\ref{act}) is now given by
\begin{equation}\label{gph}
\gamma_j. \ph_n=(-1)^{\eps_j}\ph_n,\ j=1,2,3, \quad \gamma_4 .\ph_n=(-1)^{n+\eps_4} \ph_n
\end{equation}
and
\begin{equation}
\begin{array}{ll}
\gamma_1. \phi_{l,n}^k= (-1)^{k+ \eps_1}\phi_{l,n}^k, & \gamma_2. \phi_{k,n}^l=(-1)^{l+\eps_2}\phi_{l,n}^k,\\[1ex]
 \gamma_3 .\phi_{l,n}^k = (-1)^{\eps_3}\phi_{l,n}^k,& \gamma_4  .\phi_{l,n}^k  =(-1)^{n+ \eps_4}\phi_{l,n}^k.
\end{array}\label{gPh}
\end{equation}
Furthermore,
\begin{equation}
\begin{array}{ll}
\gamma_1.\theta_{m,n,k}= (-1)^{k +\eps_1}\theta_{m,n,k}, & \gamma_2.\theta_{m,n,k}=(-1)^{\eps_2}\theta_{m,n,k+rm},\\[1ex]
\gamma_3 .\theta_{m,n,k}= (-1)^{m + \eps_3}\theta_{m,n,k},& \gamma_4 . \theta_{m,n,k}=(-1)^{n+ \eps_4}\theta_{m,n,k}.
\end{array}\label{gth}
\end{equation}
Indeed, the formula for $\gamma_2.\theta_{m,n,k}$ follows from the following observation. For $k\in\ZZ_{r'|m|}$, we have
\begin{eqnarray*} 
\lefteqn{(\gamma_2^* \theta_{m,n,k})(x,y,z,t) }\\
&=&e^{2\pi ir'm(z+\frac{\mu}{8}-\frac{\sqrt\nu}2y)}e^{int/{\kappa'}}\sum_{j\in \frac  k{r'm}+\Z}
e^{-\pi r' |m|(x-\frac{\sqrt\nu}2+j\sqrt\nu)^2}e^{\pi i r'jm(j\mu+2\sqrt\nu y-\mu)}\\
&=&e^{2\pi ir'mz}e^{int/{\kappa'}}\sum_{j\in \frac  k{r'm}+\Z}
e^{-\pi r' |m|(x+(j-\frac12)\sqrt\nu)^2}e^{\pi ir'm(\frac{\mu}{4}-\sqrt\nu y)}e^{\pi i r'jm(j\mu+2\sqrt\nu y-\mu)}\\
&=&e^{2\pi ir'mz}e^{int/{\kappa'}}\sum_{j\in \frac  k{r'm}+\Z}
e^{-\pi r' |m|(x+(j-\frac12)\sqrt\nu)^2}e^{\pi i r'(j-\frac12)m((j-\frac12)\mu+2\sqrt\nu y)}\\
&=& e^{2\pi ir'mz}e^{int/\kappa'}\sum_{j\in \frac  k{r'm}+\Z+ \frac 1{2}}
e^{-\pi r' |m|(x+j\sqrt\nu)^2}e^{\pi i r'jm(j\mu+2\sqrt\nu y)}\\
&=& \theta_{m,n, \tilde k} 
\end{eqnarray*}
with $\tilde k := k+ \frac 12 r'|m| \mod r'|m|$. 

Consequently, $L^2(\Gamma'\backslash \Osc_1^M)^\Gamma$ decomposes into $(\cH_0')^{\Gamma}\oplus(\cH_1')^{\Gamma}$, where the subrepresentations $(\cH_0')^{\Gamma}$ and $(\cH_1')^{\Gamma}$ are given by
\begin{eqnarray*}
(\cH_0')^{\Gamma}&\cong& \bigoplus_{n\in {\Bbb Z}}\bm _{\cC}(n) \cC_{\frac n{2\pi \kappa'}} \oplus \bigoplus_{a\in{\Bbb A}(\mu,\nu)}
 \bigoplus_{K=0}^{\kappa'-1} \bm_{\cS}(a,K)  {\cal S}_{a}^{K/\kappa'},\\
(\cH_1')^{\Gamma}&\cong&\bigoplus_{m\in \Z_{\not=0}} \bigoplus_{n\in\Z}\bm_\cF(m,n) \cF_{r'm,\,\frac n{2\pi\kappa'}},
 \end{eqnarray*}
where $\bm _{\cC}(n)$, $\bm_{\cS}(a,K)$ and $\bm_\cF(m,n)$ are given by (\ref{mcn}), (\ref{msa}) and (\ref{mfmn}), respectively. Indeed, the formula for $\bm _{\cC}(n)$ follows from (\ref{gph}) and the one for $\bm_{\cS}(a,K)$ from (\ref{gPh}). Finally, $\bm_\cF(m,n)$ is obtained from (\ref{gth}), where we used that $\eps_3=0$ if $r$ is odd.

In order to obtain $L^2(\Sigma(X))$, we have to pull back $(\cH_0')^{\Gamma}\oplus (\cH_1')^{\Gamma}$ by $F_S$ and to tensor the result by $\Delta$, i.e., to multiply by 4. By (\ref{Fiso}) we have
\[F_S^*\,\cC_{\frac n{2\pi \kappa'}}=\cC_{\frac n{2\pi \kappa'}},\ F_S^*\, {\cal S}_{a}^{K/\kappa'}={\cal S}_{a/2}^{K/\kappa'},\ F_S^*\,\cF_{r'm,\,\frac n{2\pi\kappa'}}=\cF_{\frac{r'm}4,\,\frac n{2\pi\kappa'}} . \]
Finally, we replace $r'$ by $2r$ and $\kappa'$ by $2\kappa$ and obtain the assertion.
\qed

\begin{re}{\rm Actually, the proof shows more than we claimed in the theorem. It gives an 
 explicit decomposition of the representation and not only an equivalence. }
\end{re}
\begin{ex}
Let us consider $X=L_r(\kappa,\mu,\nu)\backslash \Osc_1$ endowed with the trivial spin structure, i.e., $\eps_1=\dots=\eps_4=0$. Then $L^2(\Sigma(X))\cong L^2(X)\otimes \CC^4$. Indeed, Theorem~\ref{theo} implies
\begin{eqnarray*}
\bm_\cC(n)&=& \left\{\begin{array}{ll} 1, & \mbox{if $n$ is even,}  \\
0, & \mbox{else,}
\end{array}\right. \\
\bm_{\cS}(a,K)&=& \left\{\begin{array}{ll}\alpha(\mu, \nu, a/2), & \mbox{if $K$ is even,}  \\
0, & \mbox{else,}
\end{array}\right.\\
\bm_\cF(m,n)&=&\left\{\begin{array}{ll} \frac{r|m|}2, & \mbox{if $m$ and $n$ are even}, \\
0, & \mbox{else,}
\end{array}\right.
\end{eqnarray*}
which coincides with the known formulas for the decomposition of $L^2(L_r(\kappa,\mu,\nu)\backslash \Osc_1)$, see {\rm\cite{FK}}.
\end{ex}
\begin{ex}
In the case, where $r$ is even and $\eps_3=1$, we have $\bm_\cC(n)=0$ for all $n\in\ZZ$, $\bm_{\cS}(a,K)=0$ for all $a$ and $K$, and 
\[\bm_\cF(m,n)=\left\{\begin{array}{ll} \frac{r|m|}2, & \mbox{if $m$ is odd and } \eps_4=n\in\ZZ_2, \\
0, & \mbox{else.}
\end{array}\right.\]
\end{ex} 
\section{The spectrum of the cubic Dirac operator}\label{S6}
In this section, we compute the spectrum of (the closure of) the cubic Dirac operator $D^{1/3}$ on $X=L\backslash \Osc_1$ for any basic lattice $L$ of $\Osc_1$. 
We obtain that $\spec(D^{1/3})$ consists only of the point spectrum and the continuous spectrum, which we will compute in Section \ref{S61} and Section \ref{S63}, respectively.  It will turn out  that $\spec(D^{1/3})=\CC$. This will prove Theorem~\ref{T1}. Finally we determine the spectrum of the Dirac operators $D^t$ for all other $t$.
\subsection{The spectrum of $-\Omega$} \label{S61}
The operator $(D^{1/3})^2$ acts as $-\Omega$ in the first factor of the tensor product in (\ref{equr}) and trivially on the second one. On each irreducible representation, $\Omega$ acts by a scalar.  Hence, Theorem~\ref{theo} allows us to compute the spectrum of $-\Omega$ on quotients of the oscillator group by basic lattices.

If $L_r(\kappa,\mu,\nu)$ is a basic lattice, the volume of $X=L_r(\kappa,\mu,\nu)\backslash \Osc_1$ only depends on the quotient of $\kappa$ by $r$. 
For given $r,\kappa\in\NN_{>0}$ we have 
\[\beta:=\frac{\pi r}\kappa=\frac{2\pi^2}{{\rm vol}(X)}.\]

Furthermore, we define the set
\[\cA_j(\mu,\nu):=\{\pi^2 a^2\mid a\in\RR_{>0},\  \alpha_j(\mu,\nu,a)\not=0\},\ j=0,1,2. 
\]
The following theorem describes the spectrum of $-\Omega$ in dependence of the spin structure on $X$. Recall that we consider only the case $(\eps_1,\eps_2)=(0,1)$ if $\eps_1\not=\eps_2$, see Remark~\ref{isometrie}. 
\begin{theo}\label{spec} Let $L=L_r(\kappa,\mu,\nu)$ be a basic lattice and $\eps:L\to \ZZ_2$ be a homomorphism.
Then the spectrum of $(D^{1/3})^2=-\Omega$ on the spinor bundle of $X=L\backslash \Osc_{0,1}$ corresponding to $\eps$ is given in the following table:
\begin{center}
\begin{small}
\renewcommand{\arraystretch}{2}
\begin{tabular}{|c|l|l|l|l|}
  \hline
  \diagbox{$(\eps_1,\eps_2)$}{$(\eps_3,\eps_4+\kappa)$}
                 & $(0,0)$ & $(0,1)$ & $(1,0)$ & $(1,1)$ \\
  \hline
  $(0,0)$       &  $\cA_0(\mu,\nu)\cup\, 4\beta\, {\Bbb Z}$         &    $\cA_0(\mu,\nu) \cup\, 2\beta\,{\Bbb Z}$      &    $2\beta\,{\Bbb Z}$         &   $\beta\, (2{\Bbb Z}+1)$            \\
  \hline
  $(0,1)$      & $\cA_1(\mu,\nu)\cup\, 4\beta\, {\Bbb Z}$          &      $\cA_1(\mu,\nu)\cup\, 2\beta\,{\Bbb Z}_{\not=0}$       &  $2\beta\,{\Bbb Z}$            &     $\beta\, (2{\Bbb Z}+1)$          \\
  \hline
  $(1,1)$        &  $\cA_2(\mu,\nu)\cup\, 4\beta\, {\Bbb Z}$        &        $\cA_2(\mu,\nu)\cup\, 2\beta\, {\Bbb Z}_{\not=0}$       &          $2\beta\,{\Bbb Z}$     &    $\beta\, (2{\Bbb Z}+1)$           \\
  \hline
\end{tabular}
\end{small}
\end{center}
\end{theo}
\proof The Casimir operator of a representation $\sigma$ with respect to $\langle\cdot,\cdot\rangle$ equals
\[\Omega_\sigma=(\sigma_*(X))^2+(\sigma_*(Y))^2+2(\sigma_*(Z))(\sigma_*(T)).\]
A straight forward computation yields $\Omega_{\cC_d}=0$, $\Omega_{{\cal S}_a^\tau}=-4\pi^2 a^2$,
$\Omega_{\cF_{c,d}}=-2\pi c(4\pi d+1)$ for $c>0$
 and $\Omega_{\cF_{c,d}}=-2\pi c(4\pi d-1)$ for $c<0$. Therefore, $\Omega_{\cF_{\frac{rm}2,\,\frac n{4\pi\kappa}}}=-\beta m(n+ \kappa)$ for $m>0$ and 
 $\Omega_{\cF_{\frac{rm}2,\,\frac n{4\pi\kappa}}}=-\beta m(n- \kappa)$ for $m<0$.

 A representation of the form $\cC_d$ for some $d$ appears in $L^2(\Sigma(X))$ if and only if $\eps_1=\eps_2=\eps_3=0$ by (\ref{mcn}).
A representation of the form $\cS_{a/2}^{K/(2\kappa)}$  appears if and only if $\bm_{\cS}(a,K)\not=0$. Furthermore, $\bm_{\cS}(a,K)\not=0$ for some $K$ if and only if $\eps_3 = 0$ and
$$\begin{array}{ll}\alpha_0(\mu, \nu,a)\not=0, &\mbox{ if } (\eps_1,\eps_2) =(0,0),\\
\alpha_1(\mu, \nu,a)\not=0& \mbox{ if }(\eps_1,\eps_2) = (0,1),\\
\alpha_2(\mu, \nu,a)\not=0,& \mbox{ if }(\eps_1, \eps_2) = (1,1).
\end{array}
$$
Suppose that $(\eps_3,\eps_4+\kappa)=(0,0)$. Then $\bm_\cF(m,n)\not=0$ if and only if $m\not=0$ is even and $n+\kappa$ is also even. Thus the representations of the form $\cF_{c,d}$ contribute the set $4\beta \ZZ$  to the spectrum of $-\Omega$. Next suppose that $(\eps_3,\eps_4+\kappa)=(0,1)$. Then $\bm_\cF(m,n)\not=0$ if and only if $m\not=0$ is even and $n+\kappa$ is odd. In this case we get the contribution $2\beta \ZZ_{\not=0}$. Now assume that $(\eps_3,\eps_4+\kappa)=(1,0)$. Then $\bm_\cF(m,n)\not=0$ if and only if $m$ is odd and $n+\kappa$ is even. This gives the contribution $2\beta\ZZ$. Finally assume that $(\eps_3,\eps_4+\kappa)=(1,1)$. Then $\bm_\cF(m,n)\not=0$ if and only if $m$ and $n+\kappa$ are odd. In this case the representations of the form $\cF_{c,d}$ contribute the set $\beta \cdot(2\ZZ+1)$.
\qed
\begin{ex}{\rm
As an example we consider the case  
 $\nu=1,\mu=0$ and use the abbreviated notation 
$\|\cdot\|:=\|\cdot\|_{\mu,\nu}$, i.e. $\|(k,l)\|^2=k^2+l^2$ ,
$\alpha(a):=\alpha(0,1,a)$ and ${\A}:=\A(0,1)$.

Then $\A=\{a\in\RR_{>0}\mid \alpha(a)\not=0\}$ 
contains exactly those $a$ for which $a^2$ is an integer and the prime factors $q\equiv 3$ mod 4 of $a^2$ appear in even powers (Two Square Theorem), see for example \cite{H}, Theorem 366. We consider the decomposition $\A=\A_0\cup\A_1\cup\A_2$, where  
\begin{eqnarray*}
{\A}_0&:=&\{a\in \A\mid a^2\equiv 0\ {\rm mod }\ 4\}=\{\sqrt 4, \sqrt 8, \sqrt{16}, \dots\},\\
{\A}_1&:=&\{a\in \A\mid a^2\equiv 1\ {\rm mod }\ 4\}=\{\sqrt 5, \sqrt 9, \sqrt{13}, \dots\},\\
{\A}_2&:=&\{a\in \A\mid a^2\equiv 2\ {\rm mod }\ 4\}=\{\sqrt 2, \sqrt{10}, \sqrt{18}, \dots\}.
\end{eqnarray*}
Then
 $\cA_j =\{\pi^2 a^2\mid a\in\A_j\}$, $j=0,1,2$. Notice that in this example the spectrum of the spin structures  $(\eps_1,\eps_2)=(0,1)$ and $(\eps_1,\eps_2)=(1,0)$ coincide. 
 }
 \end{ex}   
\begin{pr} Let $L$ be a basic lattice and
consider a fixed spin structure corresponding to a homomorphism $\eps:L\to \ZZ_2$. The spectrum of $(D^{1/3})^2=-\Omega$ on  $X=L\backslash \Osc_1$ is symmetric if and only if  $\eps$ restricted to $L\cap Z(H)$ is non-trivial. If $\spec(-\Omega)$ is symmetric and contains $0$, then it equals $\frac{4\pi^2}{{\rm vol}(X)}\cdot\ZZ$. If it is symmetric but does not contain $0$, then it is equal to $\frac{2\pi^2}{{\rm vol}(X)}\cdot(2\ZZ+1)$.
\end{pr}
\proof Up to some inner automorphism, $L$ is equal to some $L_r(\kappa, \mu,\nu)$ for some $\kappa$, $r$, $\mu, \nu$, which we do not known explicitly. The spin structure on $L_r(\kappa,\mu,\nu)\backslash \Osc_1$ corresponding to the one on $L\backslash \Osc_1$ under this automorphism is given by a quadruple $(\eps_1,\ldots,\eps_4)$, which we also not know explicitly. Only $\eps_3$ can immediately determined by $\eps$. Indeed, $\eps_3=0$, if $\eps$ restricted to $L\cap Z(H)$ is trivial and $\eps_3=1$, if  $\eps$ restricted to $L\cap Z(H)$ is non-trivial. Now it follows from the table that the  assertions of the proposition are equivalent to the condition that neither $\cA_0$, $\cA_1$ nor $\cA_2$ is contained in $2\beta\,\ZZ$.

Since $\cA_0=\{\pi^2\|(k,l)\|_{\mu,\nu}^2\mid (k,l)\in\ZZ^2,  \mbox{ $k$ even, $l$ even }\}$, the set $\cA_0$ contains the elements $\pi s_j$, $j=1,2,3$ for 
\[\textstyle s_1:=\pi \|(0,2)\|_{\mu,\nu}^2,\ s_2:=\pi \|(2,0)\|_{\mu,\nu}^2,\ s_3:=\pi\|(2,2)\|_{\mu,\nu}^2.\]
The condition $\cA_0\subset 2\beta\, \ZZ$ would imply that $s_1, s_2$ and $s_3$ are rational. Since $s_1=4\pi /\nu$, we get $\nu=\pi q$ for some $q\in\QQ$. Furthermore, we have $s_3-s_2=s_1(1-2\mu)$, which implies that $\mu$ is rational. But then $s_2=4\pi(\nu +\frac  {\mu^2}\nu)=4(\pi^2q+\frac{\mu^2}q)$ would be irrational, which is a contradiction.

Similarly, $\cA_1\subset 2\beta\,\ZZ$ would imply that $\pi\|(0,1)\|_{\mu,\nu}^2$, $\pi\|(2,1)\|_{\mu,\nu}^2$, and $\pi\|(2,-1)\|_{\mu,\nu}^2$ are rational, which as above leads to a contradiction. Finally, if we assume that $\cA_2\subset 2\beta\,\ZZ$ we can use that $\pi\|(1,1)\|_{\mu,\nu}^2$, $\pi\|(1,-1)\|_{\mu,\nu}^2$, and $\pi\|(1,3)\|_{\mu,\nu}^2$ would be rational, which is also impossible.
\qed
\subsection{Point spectrum and eigenspaces of $D^{1/3}$}\label{S62}
In this subsection we compute  the point spectrum of $D^t$ on $X=L\backslash \Osc_1$ for  basic lattices $L$ of $\Osc_1$. For the rest of the paper we use the convention $\sqrt\lambda :=i\sqrt{|\lambda|}$ if $\lambda<0$. 

\begin{pr} Let $L$ be a basic lattice of $\Osc_1$ and consider $X:=L\backslash\Osc_1$. Then
    $z\in\CC$ belongs to the point spectrum of $D^{1/3}$ if and only if $z^2$ is in the spectrum of $-\Omega.$
\end{pr}
\proof
The operator $D^{1/3}$ preserves each summand of the decomposition of $L^2(\Sigma(X))$   obtained in Theorem~\ref{theo}. Moreover, each summand is an eigenspace of $(D^{1/3})^2=-\Omega$ with some eigenvalue $\lambda$. Let $\psi$ be an eigenspinor of $-\Omega$ with eigenvalue $\lambda$. Then $\psi_\pm := D^{1/3}\psi \pm \sqrt \lambda \psi$ satisfy the equation $D^{1/3}\psi_\pm=\pm \sqrt\lambda \psi_\pm$. One of these two spinors has to be non-trivial. Thus $\sqrt\lambda$ or $-\sqrt \lambda$ is an eigenvalue of $D^{1/3}$. Since on even-dimensional manifolds the point spectrum of the cubic Dirac operator is symmetric with respect to the imaginary axis, both $\pm\sqrt{\lambda}$ are eigenvalues of $D^{1/3}$. 
\qed

 We want to decompose the summands of the decomposition of $L^2(\Sigma(X))$ obtained in Theorem~\ref{theo} into generalised eigenspaces of $D^{1/3}$. Each summand belongs to an eigenspace of $-\Omega$. Let $\lambda$ be the corresponding eigenvalue. If  $\lambda\not=0$ the summand decomposes into eigenspaces of $D^{1/3}$ with eigenvalues $\pm \sqrt \lambda$. The projections to these subspaces are given by $P_\lambda^\pm:=\frac{1}{2\sqrt \lambda}(\pm D^{1/3} + \sqrt \lambda)$.
For the summands with $\lambda=0$, we will determine the kernel of $D^{1/3}$.

We will use the notation of  subsection \ref{S34} for the irreducible unitary representations of the oscillator group. Furthermore, we consider  the complete orthonormal systems of the representation spaces  introduced there. We also use the formulas for the Dirac operator  with respect to the basis $(u_1,\dots, u_4)$ of $\Delta$ introduced in subsection \ref {Diracosc}. 

Let $\cF_{c,d}\otimes\Delta$, $c>0$, be one of the summands in (\ref{H1}). Assume first that $\lambda\not=0$. Then the eigenvalue of $-\Omega$ on $\cF_{c,d}$ equals $\lambda=2\pi c(4\pi d+1)$. Using the projections $P_{\pm}$ we calculate that   the 
 subspace of  eigenvectors of $D^{1/3}$ with eigenvalue $\pm\sqrt\lambda$ is spanned by
\begin{eqnarray}
&&\eta_{0}^\pm=\sqrt{2\pi c} \,\psi_0\otimes u_1 \pm i\sqrt{2\pi d+\textstyle\frac12}\,\psi_0\otimes u_2, \label{ev1}
\\
&&\eta_{n}^\pm=\sqrt{2\pi c}\, \psi_{n} \otimes u_1 \pm i\sqrt{2\pi d+\textstyle\frac12}\, \psi_{n}\otimes u_2+\sqrt{n}\,\psi_{n-1}\otimes u_4,\quad n\ge 1, \label{ev2}\\
&&\hat \eta_{n}^\pm =\sqrt{n+1}\, \psi_{n+1}\otimes u_2 -\sqrt{2\pi c}\, \psi_n \otimes u_3 \mp i\sqrt{2\pi d+\textstyle\frac12}\, \psi_n\otimes u_4,\quad n\ge 0.\label{ev3}
\end{eqnarray}
If $\Omega=0$ on $\cF_{c,d}$ or equivalently if $4\pi d+1=0$, then $\Ker D^{1/3}= \Im D^{1/3}$ on $\cF_{c,d}$. This space is spanned by
\begin{eqnarray*}
&&\eta_{0}=\sqrt{2\pi c} \,\psi_0\otimes u_1, \\
&&\eta_{n}=\sqrt{2\pi c}\, \psi_{n} \otimes u_1 +\sqrt{n}\,\psi_{n-1}\otimes u_4,\quad n\ge 1,\\
&&\hat \eta_{n}= \sqrt{n+1}\, \psi_{n+1}\otimes u_2 -\sqrt{2\pi c}\, \psi_n \otimes u_3 ,\quad n\ge 0. 
\end{eqnarray*}
 Now let $\cF_{c,d}\otimes\Delta$, $c<0$, be one of the summands in (\ref{H1}). Then the eigenvalue of $-\Omega$ on $\cF_{c,d}$ equals $\lambda=2\pi c(4\pi d-1)$.
If $4\pi d-1>0$, the subspace of  eigenvectors of $D^{1/3}$ with eigenvalue $\pm\sqrt{\lambda}$  is spanned by 
\begin{eqnarray}
&&-\sqrt{2\pi |c|} \,\psi_0\otimes u_3 \mp  \sqrt{2\pi d-\textstyle\frac12}\,\psi_0\otimes u_4,\label{ev4}\\
&& \sqrt{2\pi |c|}\, \psi_{n} \otimes u_1 \pm \sqrt{2\pi d-\textstyle\frac12}\, \psi_{n}\otimes u_2+\sqrt{n+1}\, \psi_{n+1}\otimes u_4,\quad n\ge 0,\label{ev5} \\
&&\sqrt{n}\,\psi_{n-1}\otimes u_2  - \sqrt{2\pi |c|}\, \psi_{n} \otimes u_3 \mp \sqrt{2\pi d-\textstyle\frac12}\, \psi_{n}\otimes u_4,\quad n\ge 1.\label{ev6}
\end{eqnarray}
If $4\pi d-1<0$, the subspace of  eigenvectors of $D^{1/3}$ with eigenvalue $\pm\sqrt{\lambda}$  is spanned by 
\begin{eqnarray}
&&-\sqrt{2\pi |c|} \,\psi_0\otimes u_3 \pm i  \sqrt{\textstyle\frac12-2\pi d}\,\psi_0\otimes u_4,\label{ev7}\\
&& \sqrt{2\pi |c|}\, \psi_{n} \otimes u_1 \mp i\sqrt{\textstyle\frac12-2\pi d}\, \psi_{n}\otimes u_2+\sqrt{n+1}\, \psi_{n+1}\otimes u_4,\quad n\ge 0,\label{ev8} \\
&&\sqrt{n}\,\psi_{n-1}\otimes u_2  - \sqrt{2\pi |c|}\, \psi_{n} \otimes u_3 \pm i \sqrt{\textstyle\frac12-2\pi d}\, \psi_{n}\otimes u_4,\quad n\ge 1.\label{ev9}
\end{eqnarray}

If $\Omega =0$, then we obtain that $\Ker D^{1/3}= \Im D^{1/3}$ on $\cF_{c,d}$ is spanned by
\begin{eqnarray*}
&&\sqrt{2\pi |c|} \,\psi_0\otimes u_3 ,\\
&&\sqrt{2\pi |c|}\, \psi_n \otimes u_1 +\sqrt{n+1}\, \psi_{n+1}\otimes u_4  \, \quad n\ge 0,\\
&&\sqrt{n}\,\psi_{n-1}\otimes u_2 -\sqrt{2\pi |c|}\, \psi_{n} \otimes u_3 ,\quad n\ge 1.
\end{eqnarray*}

On summands of the form $\cS_a^\tau\otimes \Delta$, $a>0$, the operator 
$-\Omega$ has eigenvalue $4\pi^2 a^2$. Using again the projections  $P_\pm$, we see that the subspaces of eigenvectors of $D^{1/3}$ with eigenvalue $\pm 2\pi a$ are spanned by 
\begin{eqnarray}
&& \phi_n\otimes u_2\mp \phi_{n+1} \otimes u_4,\quad n\in\ZZ, \label{neu1}
\\
&&\pm \phi_n\otimes u_1+\textstyle \frac i{\sqrt2 \pi a}(n+\tau +\frac12)\, \phi_n\otimes u_2 +\phi_{n+1}\otimes u_3,\quad n\in\ZZ .\label{neu2}
\end{eqnarray}
All summands of the form $\cC_d$ are in the kernel of $D^{1/3}$.

Let $\lambda\not=0$ be an eigenvalue of $-\Omega$ on $L^2(\Sigma(X))$ and $E_\lambda$ be the corresponding eigenspace.
We have seen that the two projections $P^{\pm}_{\lambda}:=\frac{1}{2\sqrt \lambda}(\pm D^{1/3} + \sqrt \lambda):E_\lambda\rightarrow E_\lambda$ to the generalised eigenspaces of $D^{1/3}$ are non-trivial. 

\begin{pr}\label{pnc} For any basic lattice $L$ and any spin structure on $X$ there exists an eigenvalue $\lambda\not=0$ of $-\Omega$ such that the projections $P_{\lambda}^\pm$ are not continuous. 
\end{pr}
\proof By Theorem~\ref{theo}, the decomposition of $L^2(\Sigma(X))$ contains a summand of the form $\cF_{c,d}\otimes \Delta$, $c>0$, for $\lambda=2\pi c(4\pi d+1)\not=0$. By (\ref{ev2}), the projections of $\zeta_n:=2i\sqrt{2\pi d+\textstyle\frac12}\, \psi_{n}\otimes u_2$ are equal to 
\[ P^\pm_{\lambda}\zeta_n = \pm\sqrt{2\pi c}\, \psi_{n} \otimes u_1 + i\sqrt{2\pi d+\textstyle\frac12}\, \psi_{n}\otimes u_2\pm\sqrt{n}\,\psi_{n-1}\otimes u_4.
\]
Since all $\zeta_n$ have the same length, this shows that the projections $P^\pm_{\lambda}$ are unbounded. 
\qed 
\subsection{The residual and the continuous spectrum}\label{S63}
\begin{pr} \label{rscs} Let $L$ be a basic lattice of $\Osc_1$. 
 The residual spectrum of $D^{1/3}$ on $X=L\backslash \Osc_1$  is empty for every spin structure on $X$. The continuous spectrum is equal to $\CC\setminus \spec_p(D^{1/3}) $. 
\end{pr}
\proof  To see that there is no residual spectrum, we use the symmetry properties of the spectrum.  We already have seen in Section \ref{S25} that $iD^t$ is essentially selfadjoint. Hence Fact \ref{fact} applies to this operator.  We obtain that the residual spectrum  $\spec_{r}(D^t)$ is contained in the complex conjugate of  $-\spec_p(D^t)$.   
The point spectrum of $\spec_p(D^{1/3})$ contains only real and purely imaginary values since the spectrum of the Casimir operator is real.  Because  $X$ has even dimension,  $\spec_p(D^{1/3})$ is symmetric to zero. In particular,  the point spectrum is  invariant under complex conjugation. Thus $\spec_{r}(D^{1/3})\subset -\spec_{p}(D^{1/3})=\spec_{p}(D^{1/3})$. Hence the residual spectrum is empty. We remark that the assumptions of \cite[Satz 3.20]{B} are satisfied, see Section \ref{S25}. Application of item 4 of this theorem shows directly that $\spec_r(D^{1/3})$ is empty. 

Since the residual spectrum is empty, it suffices to show that the approximate spectrum of $D^{1/3}$ is equal to $\CC$.  Take $z\in \CC \setminus \spec_p ( D^{1/3})$. We  have to  show that there is a sequence of spinors $(\Phi_j)_{j \in {\Bbb N}}$ with $\| \Phi_j\| = 1$ and $\| (D^{1/3}-z )\Phi_j \| \to 0$. Note first that for any $\lambda\in \spec_p (( D^{1/3})^2)$ and any corresponding unit eigenspinor $\Psi$ of $(D^{1/3})^2$ the unit spinor $\Phi =\| (D^{1/3}+z )\Psi \|^{-1} (D^{1/3}+z )\Psi$ satisfies 
 $\| (D^{1/3}-z )\Phi\| = \| (D^{1/3}+z )\Psi \|^{-1}  \vert\lambda - z^2\vert.$ Therefore it suffices to find a sequence of unit eigenspinors $\Psi_n$ of $(D^{1/3})^2$ such that $ \| (D^{1/3}+z )(\Psi_n) \|  \to \infty$. Recall that $L^2(\Sigma(X))$ contains summands equivalent to $\cF_{c,d}\otimes \Delta$ for some $c>0,d\in\RR$ with $\lambda=2\pi c(4\pi d+ 1)\not=0$, see \ref{theo}. Therefore the existence of such a sequence follows from the unboundedness of the projections $P^\pm_{\lambda} $ in this summand. \qed
\subsection{A remark on the point spectrum of quotients by shifted lattices}\label{S64}
In subsections~\ref{S61} and \ref{S62} we proved that for basic lattices in $\Osc_1$ the point spectrum of the cubic Dirac operator on the quotient space is discrete.
Here we will show that there also exist (non-basic) lattices such that the point spectrum is not discrete. A similar statement for the wave operator has been proven in $\cite{FK}$. 
\begin{ex}{\rm
Let $L=L_r(\kappa,\mu,\nu)$ be a basic lattice. We choose a real number $u$ such that $\tilde u:=2\pi \kappa r u$ is irrational and consider $L':=T_u(L)$, where $T_u$ is the automorphism defined in (\ref{DTu}). We consider a spin structure for which $\eps_3=0$ and $\eps_4=\kappa$.   By (\ref{pull2}), the right regular representation $L^2(\Sigma(L'\backslash \Osc_1))$ is equivalent to $(T_u^{-1})^*( L^2(\Sigma(L\backslash\Osc_1))$. If $m\in\NN_{>0}$ and $n\in\ZZ$ satisfy $(m,n)=(\eps_3,\eps_4)=(0,\kappa)\in\ZZ_2^2$, then $\cF_{\frac{rm}2,\frac n{4\pi\kappa}}$ appears as a summand in $L^2(\Sigma(L\backslash\Osc_1)))$, see Theorem~\ref{theo}. 
Hence, in this case, $(T_u^{-1})^*(\cF_{\frac{rm}2,\frac n{4\pi\kappa}})=\cF_{\frac{rm}2,\frac n{4\pi\kappa}-u\frac{rm}2}$ appears as a summand in $L^2(\Sigma(L'\backslash \Osc_1))$. In particular, 
the point spectrum of $\Omega$ on $\Sigma(L'\backslash \Osc_1)$ contains the eigenvalue of $\Omega_{\cF_{\frac{rm}2,\frac n{4\pi\kappa}-u\frac{rm}2}}$, which equals $-\beta m(n+\kappa-\tilde u m)$. Therefore $\spec_p(\Omega)$ contains the set 
\begin{eqnarray*}
\cB&:=&\{-\beta m(n+\kappa-\tilde u m)\mid m\in 2\NN_{>0},\ n+\kappa\in 2\ZZ\}\\
&=&\{-4\beta m'(n'-\tilde u m')\mid m'\in \NN_{>0},\ n'\in\ZZ\}
\end{eqnarray*}
By Dirichlet's approximation theorem, for every $N\in\NN$, there exists a pair $(m',n')\in \NN\times\ZZ$ with $0<m'\le N$  such that 
\[\left| n'-\tilde u m'\right|<\frac 1{N}.\]  This implies that the set 
$\{(n'-\tilde u m')m'\,|\,m'\in\NN_{>0},\,n'\in\ZZ\}$ contains infinitely many numbers in $[-1,1]$ since $\tilde u$ is irrational. Hence, $\cB$ and therefore also $\spec_p(\Omega)$ contains an accumulation point in $\RR$. 
}\end{ex}

\subsection{Application: the spectrum of $D^t$}
\label{S65}
In this subsection we consider the Dirac operator $D^t$ for arbitrary $t\in \RR$.
\begin{pr} Let $L$ be a basic lattice of $\Osc_1$ and consider $X:=L\backslash\Osc_1$. 
    The Dirac operator respects the decomposition of $L^2(\Sigma(X))$ into irreducible summands given in Theorem \ref{theo}. The eigenvalues of $D^t$ on the various summands of type $\cF_{c,d}$, $\cS_a$, and $\cC_d$ are given in the following table:
\[
\mbox{   
\begin{small}
\renewcommand{\arraystretch}{2}
\begin{tabular}{|c||c|c|}
  \hline
  rep & eigenvalues of $D^{1/3}$ &eigenvalues of $D^t$\\
  \hline \hline
$\cF_{c,d}$, $c>0$  & $\pm\sqrt\lambda,\ \lambda=2\pi c(4\pi d+1)$
  & \\
  \cline{1-2}
   $\cF_{c,d}$, $c<0$  & $\pm\sqrt\lambda,\ \lambda=2\pi c(4\pi d-1)$  &\raisebox{3ex}[-3ex]{$\pm (\lambda \pm 2\pi c (3t-1))^\frac12$ }  \\
  \hline
    $\cS_{a}$&$\pm 2\pi a$
  &  $\pm 2 \pi a$\\
  \hline
 $\cC_{d}$&$0$&$0$ \\
  \hline
\end{tabular}
\end{small}
}
\]
The spectrum of $D^t$ is equal to $\CC$ and the residual spectrum is empty. 
\end{pr}
\proof 
According to (\ref{xyz}), we have
\begin{equation}\label{Dt}
 \textstyle D^t = D^{1/3} +   \frac 12 i(3t-1)\,  \Id \otimes  \left(\begin{array}{cc} A&0 \\0 & -A\end{array}\right).
\end{equation}
Therefore it respects the decomposition of $L^2(\Sigma(X))$. Let us first consider the restriction of $D^t$ to summands of type $\cF_{c,d}$ for $c>0$, $4\pi d+1\not=0$.  Recall that we obtained eigenvectors $\eta_n^\pm$, $\hat \eta_n^\pm$ of $D^{1/3}$ given by (\ref{ev1})--(\ref{ev3}).
The subspaces $V_n:=\Span\{\eta_n^+, \eta_n^-\}$ and  $\hat V_n:=\Span\{\hat \eta_n^+, \hat \eta_n^-\}$, $n\ge0$, are invariant under $D^t$ and their sum spans $\cF_{c,d}$. 
On these subspaces, $D^t$ is given by the matrix
\begin{equation} 
 D^t =  \left(\begin{array}{cc} \sqrt{\lambda} +b&b \\-b &-( \sqrt{\lambda} +b)\end{array}\right)
\end{equation}
with $b = \pm (3t-1) \sqrt{\frac{\pi c}{8 \pi d + 2}}$, where the plus sign appears on $V_n$ and the minus on $\hat V_n$. Thus the eigenvalues given in the table can be calculated. The formulas for $c<0$ and $4\pi d-1\not=0$ follow  in the same way using the eigenvectors given in (\ref{ev4})--(\ref{ev9}).

In the case $c>0$, $4\pi d+1=0$, $\cF_{c,d}$ also splits into two-dimensional invariant subspaces of $D^t$. Indeed, $D^t$ leaves invariant $\Span\{\eta_n, \zeta_n\}$ and $\Span\{\hat \eta_n, \hat \zeta_n\}$  with $\zeta_n = \psi_n \otimes u_2$ and $\hat \zeta_n = \psi_n \otimes u_4$ for $n \geq 0$. It is given by
\begin{equation}\label{0}
D^t = i \sqrt{ \pi c} \left(\begin{array}{cc} 0& \mp 2 \\ 3t-1 &0\end{array}\right).
\end{equation} 
with the minus sign on $\Span\{\eta_n, \zeta_n\}$ and the plus sign on  $\Span\{\hat \eta_n, \hat \zeta_n\}.$
For $c<0$, $4\pi d-1=0$, we obtain also a splitting into two-dimensional invariant subspaces on which $D^t$ is again given by (\ref{0}).  

On summands of the form $\cS_a^\tau\otimes \Delta$, $a>0$, the vectors in (\ref{neu1}) and (\ref{neu2}) span an invariant subspace for every $n\in\ZZ$ with $D^t$ given by 
\begin{equation}
 D^t = \pm \left(\begin{array}{cc} 2 \pi a&s \\ 0  & 2 \pi a\end{array}\right),
\end{equation}
where $ s = \frac{ i  \sqrt{2}}{2} (3t-1). $ 

On each summand of type $\cC_d$ the cubic Dirac operator acts trivially. Hence  equation~(\ref{Dt}) implies that $D^t$ is nilpotent on this summand.  Thus its only eigenvalue is~0.

There is no residual  spectrum for the same reasons as for $D^{1/3}$ because the point spectrum again contains only real and  purely imaginary values.

The continuous spectrum is equal to $\CC\setminus \spec_p(D^{t})$. To see this we use the same argument as for $D^{1/3}$. We can choose $c>0$ and $d\in\RR$ with $\lambda_t:=2\pi c(4\pi d +3t)\not=0$ such that $\cF_{c,d}\otimes\Delta$ appears in the decomposition of $L^2(\Sigma(X))$. As in the proof of Proposition \ref{rscs} we have only to show that 
the  two projections 
$ P^{\pm}_t:=\frac{1}{2\sqrt{\lambda _{t}}}(\pm D^{t} + \sqrt{\lambda_t})$ from $\cF_{c,d}\otimes\Delta$ to the generalised 
eigenspaces 
 of $D^{t}$ with eigenvalues $\pm\sqrt{\lambda_t}$ are not continuous. However, 
$\sqrt{\lambda _{t}}P_t^\pm$ differs from $\sqrt{\lambda_{1/3} }P _{1/3}^\pm$ only by a bounded operator and the assertion follows from Proposition~\ref{pnc}. \qed

\end{document}